\documentclass[11pt]{amsart}
\usepackage{amssymb,amscd}
\usepackage[all,cmtip]{xy}
\usepackage{epic}
\usepackage{url}
\numberwithin{equation}{section}
\newtheorem{theorem}{Theorem}[section]
\newtheorem{lemma}[theorem]{Lemma}
\newtheorem{proposition}[theorem]{Proposition}
\newtheorem{corollary}[theorem]{Corollary}

\newtheorem{definition-lemma}[theorem]{Definition-Lemma}
\theoremstyle{definition}

\newtheorem{example}[theorem]{Example}

\theoremstyle{remark}
\newtheorem{remark}[theorem]{Remark}

\def\FF{{\mathbb F}}

\def\PP{{\mathbb P}}
\def\QQ{{\mathbb Q}}

\def\ZZ{{\mathbb Z}}

\def\cO{{\mathcal O}}

\begin{document}
\title[Optimal curves of low genus over finite fields]
{Optimal curves of low genus over finite fields}

\author{Alexey Zaytsev}

\address{ Immanuel Kant Baltic Federal University, Nevsky 14a, Kaliningrad, Russia}

\email{alzaytsev@kantiana.ru}

\begin{abstract}
The Hasse-Weil-Serre bound is improved for curves of low genera over 
finite fields with  discriminant in $\{-3,-4,-7,-8,-11,-19 \}$ by
studying  optimal curves.
\end{abstract}
\maketitle
\begin{section}{Introduction.}
In this paper we give some improvements of the Hasse-Weil-Serre upper and
lower bound for the number of rational points on a curve over a finite field
for curves of low genus over finite fields of certain discriminants.

Here by a curve over a finite field $\FF_q$ we mean an absolutely 
irreducible nonsingular projective algebraic variety of dimension $1$ 
over $\FF_{q}$ and by the discriminant $d(\FF_{q})$ of a finite 
field  $\FF_q$ we mean  the integer $m^2-4q$, where $m=[2\sqrt{q}]$. 
Throughout the paper we denote the set $ \{-3,-4,-7,-8,-11,-19 \}$ by $D$.

For the number of rational points on a curve $C$ of genus $g$ 
defined over a finite field $\FF_q$ we have the well-known bound of 
Hasse-Weil as improved by J-P.\ Serre
$$
|\#C(\FF_q)-q-1| \le g[2\sqrt{q}],
$$
which now is called the Hasse-Weil-Serre bound.

This bound is interesting in its own right, but also in view
of applications in coding theory, as for example for the geometric
codes introduced by V.\ D.\ Goppa in 1980.

A curve $C$ of genus $g$ over a finite field $\FF_q$ is called a
maximal (resp.\ minimal) optimal curve if its number of rational points
attains the upper (resp.\ lower) Hasse-Weil-Serre $q+1\pm g[2\sqrt{q}]$.

We study optimal curves of low genera over fields with discriminant
$d(\FF_q) \in D$. It turns out that such curves
are ordinary. This allows us to use the canonical lifts of the Jacobian
and the equivalence of categories between the category of ordinary principally
polarized abelian varieties over $\FF_q$ and the category 
of certain unimodular irreducible
hermitian modules over the ring of integers $\cO_K$ of an imaginary quadratic
field $K$ of discriminant equal to $d(\FF_q)$. The use of such hermitian modules
to study curves over finite fields was initiated by J-P.\ Serre \cite{Serre2}, \cite{Serre3} and 
continued by K.\ Lauter \cite{Lauter1}. For the discriminants we consider here
the class number of $\cO_K$ is $1$ and there is a classification of
hermitian modules. This together with results on defect $1$ curves
leads to our improvement.

\begin{theorem} Let $C$ be a curve of genus $g$ over a finite field $\FF_q$
of characteristic $p$. 
Then we have that
$$
|\#C(\FF_q)-q-1| \le g[2\sqrt{q}]-2,
$$
if the conditions on $q$ and $g$ given by a line in the following 
table hold:
\begin{center}
\begin{tabular}{|c|c|r|l|c|}
\hline
$d(\FF_q)$ & $q$ & $g$  \\
\hline
$-3$ & $q\neq 3$ & $3\leq g \leq 10$ \\
$-4$ & $q\neq 2$ & $3\leq g \leq 10$ \\
$-7$ &  & $4 \leq g \leq 8$ \\
$-8$ & $p\neq 3$ & $3\leq g \leq 7$ \\
$-11$ & $p\neq 3,\,  q < 10^4 $ & $g=4$  \\
$-11$ & $p>5$ & $g=5$  \\
$-19$ & $q< 10^3$ and $q \equiv 1 \, (\bmod\, 5)$  & $g=4$ \\
\hline
\end{tabular}
\end{center}
\end{theorem}
\par
Our method uses the explicit classification of hermitian lattices  by A.~Schiemann \cite{Sch}.
We also use some information on generators for the automorphism groups of such lattices of 
the dimension $4$ and $5$ over an imaginary quadratic extension $K$ of $\QQ$, with discriminant ${d}(K)=-19$
provided to us by R.\ Schulze-Pillot \cite{Schulze-letter}.
\subsection*{Acknowledgments} The author would like to thank  G.\ van der Geer for his generous help and 
suggestions.  The author is also grateful to  R.\ Schulze-Pillot 
for providing the generators of the automorphism groups of unimodular hermitian modules, which were obtained with  use of a computer program, developed by A.\ Schiemann. The research of the author has been financed by the Netherlands Organisation for the Advancement of Scientific Research (NWO).
\end{section}
\begin{section}{An equivalence of categories.}
Let $C$ be a maximal or minimal curve of genus $g$ over $\FF_{q}$. Then
by a corollary of the Hasse-Weil-Serre bound (\cite{Shabat} Theorem 5) and by the Honda-Tate theorem \cite{Waterhouse}
the Jacobian ${\rm  Jac}(C)$ is isogenous to $E^g$, where $E$ is a maximal or a minimal elliptic curve, respectively. 
If the elliptic curve $E$ is ordinary, the Jacobian  ${\rm Jac}(C)$ is an ordinary 
principally polarized abelian variety and we can apply an equivalence of categories 
between the category of ordinary abelian varieties over a finite field  
and the category of Deligne modules; see \cite{Deligne}, 
\cite{Howe}, or \cite{Serre} for a more concrete description 
of this equivalence. The idea is that the canonical lift of our Jacobian is a
complex abelian variety with a certain endomorphism ring
and in this way we reduce to a classification problem of hermitian modules.
\par
In order to use the classification of 
unimodular irreducible hermitian modules as given by A.~Schiemann in \cite{Sch},
we have to work  with maximal orders. According to a result in \cite{Deligne} the 
isomorphism classes of the abelian varieties 
in the isogeny class [$E^g$] correspond one-to-one 
to the isomorphism classes of $R$-modules that can be embedded 
as a lattice in the $K$-vector space $K^g$, where $R={\rm  End}(E)$ and $K={\rm  Quot}(R)$.
Since our optimal elliptic curve is defined over a finite field $\FF_{q}$ with the discriminant
$d \in D$ we have that $R=\cO_K$ is the ring of integers of $\QQ(\sqrt{d})$,
with class number $1$ and there is only one isomorphism class of such an $R$-modules. 
Hence there is only one isomorphism class of abelian varieties 
in the isogeny class [$E^g$], so ${\rm  Jac}(C)\cong E^g$ as abelian varieties over $\FF_q$.
According to \cite{Serre} the principal polarizations of 
${\rm  Jac}(C)$ correspond to unimodular hermitian 
$\cO_{K}$-forms $h: \cO_{K}^g \times \cO_{K}^g \rightarrow \cO_{K}$. 
Moreover, we have  ${\rm  Aut}_{\FF_{q}}({\rm  Jac}(C))\cong {\rm  Aut}(\cO_{K}^g,h)$, where ${\rm  Aut}(\cO_{K}^g,h)$ is the automorphism group of the hermitian module $\cO_{K}^g$ in $(K^g,h)$. 
This allows us to use the classification of unimodular irreducible hermitian forms with discriminant $d(K)$ and their automorphism groups.
\par
We need conditions on the finite fields $\FF_q$ such that there can exist an optimal
ordinary elliptic curve over $\FF_q$.
\par
\begin{proposition}Let $E$ be an optimal elliptic curve over a finite field $\FF_{q}$ with discriminant
$d(\FF_q) \in D$. Then $E$ is ordinary if $q \neq 2,\; 3$. 
\end{proposition}
\begin{proof}
Let $E$ be an optimal elliptic curve over the finite field $\FF_{q}$ and  $d(\FF_q) \in D$.
If $E$ is supersingular, then ${\rm  char}(\FF_q)$ has to be a divisor of ${d}(\FF_q)$, 
since ${\rm  char}(\FF_{q})$ divides $m$ and hence ${\rm  char}(\FF_{q})$ divides ${d}(\FF_{q})=m^2-4 \cdot q$.
Now we explore each discriminant separately.
\par
If $E$ is supersingular and ${d}(\FF_q)=-3$ then ${\rm  char}(\FF_q)=3$ and $m^2=3(4\cdot 3^{n-1}-1)$. 
Then $3$ divides $4\cdot 3^{n-1}-1$, hence $n=1$. 
\par
If $d(\FF_{q})=-4$ then the optimal elliptic curve $E$ over $\FF_q$ is not supersingular.
Otherwise, by $m^2=4q-4$ it follows that ${\rm char}(\FF_{q})=2$ and $(m/2)^2=2^{n-1}-1$ 
with $m/2 \in \ZZ$, hence $(m/2)^2 = 2^{n-1}-1 \equiv 3 \; (\mbox{{\rm mod}}\; 4)$, but $3$ is not a square in $\ZZ/4\ZZ$. 
\par
The optimal elliptic curve $E$ is ordinary if  the discriminant of the finite field $\FF_q$ is in $\{-7, -11, -19\}$. 
In fact, if $E$ is  supersingular then it is defined over $\FF_{{|d|}^n}$, where $d:=d(\FF_q)$, and hence  $m^2=|d|\cdot(4(|d|)^{n-1}-1)$ but this is impossible, since $d$ is prime.

For  discriminant $d(\FF_{q})=-8$ we get a contradiction in a similar way.
\end{proof}
We now have the following corollary.
\begin{corollary}
If $C$ is an optimal curve over a finite field $\FF_{q}$ with discriminant $d(\FF_q) \in D$ then ${\rm  Jac}(C)$ is an ordinary abelian variety if $q \neq 2,\; 3$.
\end{corollary}

Since our curve embeds in its Jacobian ${\rm  Jac}(C)\cong E^q$
the projections on the factors $E$ define morphisms $C \to E$ and
we want to calculate  the degree of these maps. Here we use
the fact that the lattice of the canonical lift is a free $\cO_K$-module, since
the class number is $1$.
\begin{proposition}\label{proj}
Let $C$ be an optimal curve over $\FF_q$ and 
fix an isomorphism ${ \rm  Jac}(C)\cong E^g$ such that the theta 
divisor corresponds to the hermitian form $(h_{ij})$ on $\cO_{K}^g$
on the canonical lift of ${\rm  Jac}(C)$. 
Then the degree of the $k$-th projection 
\begin{equation*}
f_{k}:C\hookrightarrow {\rm  Jac}(C)\cong E^g \stackrel{pr_k}{\longrightarrow} E
\end{equation*}
equals $\det( h_{ij})_{i,j \neq k}$.
\end{proposition}
\begin{proof}
We denote the abelian variety $E^g$ by
 $E_{1}\times \ldots \times E_{g}$, where $E_{i}=E$,
and consider the first projection. 
The degree of the map $f_1$ equals the intersection 
number $[C]\cdot[E_2 \times \ldots \times E_{g}]$. 
The cohomology class $[C]$ of $C$ in an appropriate cohomology theory
is $[\Theta^{g-1}/(g-1)!]$. Recall that
if $L$ is a line bundle on an abelian variety $A$ of dimension $g$ 
then by the Riemann-Roch theorem one has $(L^g/g!)^2=\deg(\varphi_{L})$,
and  $\deg(\varphi_{L})=\det(r_{ij})^2$, where the matrix
$(r_{ij})$ gives the hermitian form corresponding to the first Chern class
of the line bundle $L$.
Since the hermitian form $(h_{ij})_{i,j \neq 1}$ corresponds 
to the line bundle $\Theta|_{E_2 \times \ldots \times E_{g}}$
on the abelian variety $E_2 \times \ldots \times E_{g}$ the degree of 
$f_1$ is given by  
$$
[C]\cdot[E_2 \times \ldots \times E_{g}]=
\frac{1}{(g-1)!}(\Theta|_{E_2 \times \ldots \times E_{g}})^{g-1}=
\det((h_{ij})_{i,j \neq 1}).
$$
\end{proof}
\end{section}
\begin{section}{Optimal curves of low genus over finite fields with discriminants $-3, -4, -7$ and $-8$.}
In this section we prove the non-existence of certain optimal curves using properties of the automorphism
groups of abelian varieties over the finite fields $\FF_q$ with the discriminant 
${d}(\FF_q) \in \{-3,\; -4,\; -7,\; -8\}$.

Throughout this section we denote by $C$  an optimal curve of genus $g$ over $\FF_{q}$ and by $h$ an unimodular irreducible hermitian form which corresponds to the polarization of ${\rm Jac}(C)$ under the equivalence of Section~2.

Using Torelli's theorem (cf.\ \cite{Andreotti}) and an upper bound on the number automorphisms of a curve over a finite field (cf.\ \cite{Roquette}) we obtain an upper bound on the number of automorphisms of irreducible unimodular hermitian modules $(\cO_{K}^g,h)$.

If a curve $C$ defined over a finite field $\FF_q$ is not exceptional 
(see Subsection~3.1) and $p>g(C)+1$ then we have
$\# {\rm  Aut}_{\FF_{q}}(C)\le 84(g(C)-1)$. Since $\# {\rm   Aut}_{\FF_{q}}({\rm  Jac}(C))=\# {\rm  Aut}_{\FF_{q}}(C)$
or $2\; \# {\rm  Aut}_{\FF_{q}}(C)$ depending on whether the curve is hyperelliptic or not, we get
\begin{displaymath}
\# {\rm Aut}(\cO_{K}^g,h) \le
\left\{ \begin{array}{ll}
84(g(C)-1)& \textrm{if}\; C\; \textrm{ is hyperelliptic,}\\
168(g(C)-1) & \textrm{otherwise.}\\
\end{array}\right.
\end{displaymath}
\par
We shall use this bound and other properties of the automorphism group of a unimodular hermitian module 
to restrict the possibilities for optimal curves of low genus. 

To begin with, we treat the case of exceptional curves and their properties. 
\begin{subsection}{Optimal Exceptional Curves.}
Recall that by an exceptional curve we mean a hyperelliptic curve over the finite field $\FF_{q}$ 
given by equation of the form $y^2=x^p-x$,
with $p={\rm  char}(\FF_{q})$. An exceptional curve $C$ has genus $g=(p-1)/2$ and 
$\# {\rm  Aut}_{\bar{\FF}_q}(C)=2p(p^2-1)$. Optimal exceptional curves over $\FF_{q}$ of genus $ 3 \le g \le 10$
can occur only if ${\rm  char}(\FF_{q}) \in \{7,\; 11,\; 13,\;17,\;19 \}$.
Now we list the possibilities for $q,\, g$ and $\#{\rm  Aut}_{\bar{\FF}_{q}}(C)$.
\renewcommand{\arraystretch}{1.25}
\begin{center}
\begin{tabular}{|c|c|l|} 
\hline
${q}$&$g(C)$&$\#{\rm  Aut}_{\bar{\FF}_{q}}(C)$\\
\hline 
${7^n}$&$3$&$2^5\cdot3\cdot7$\\
${11^n}$&$5$&$2^2\cdot3^2\cdot5\cdot11$\\
${13^n}$&$6$&$2^4\cdot3\cdot7 \cdot13$\\
${17^n}$&$8$&$2^6\cdot3^2\cdot17$\\
${19^n}$&$9$&$2^4\cdot3^2\cdot5 \cdot 19$\\
\hline
\end{tabular}
\end{center}
The order of the automorphism group ${\rm  Aut}_{\FF_{q}}(C)$ 
then equals $\#{\rm  Aut}(\cO_{K}^{g},h)$ for some hermitian module $(\cO_{K}^{g},h)$ and it also has to divide 
$\#{\rm  Aut}_{\bar{\FF}_{q}}(C)$.
However from Tables~1,2,3,4 (see here after) for the orders of the automorphism groups of such modules one sees that
the order of the automorphism group of any unimodular irreducible module of dimension $g(C)$ does not divide $\#{\rm  Aut}_{\bar{\FF}_{q}}(C)$ (under the restriction of each subsection), and hence an exceptional curve cannot be an optimal curve
in the case we consider here.

There is another argument to exclude exceptional curves which is base on $p$-rank.
The Deuring-Shavarevich  formula (see \cite{Nakajima} ) yields that an exceptional curve is never ordinary, 
and hence it should be excluded from the consideration. 
\end{subsection}
\par
In the remainder of this section we may and shall assume that our optimal curves are not exceptional curves. 
\begin{subsection}{Optimal curves over finite fields with discriminant $-3$.}
Here we consider optimal curves of genus $g \le 10$ over finite fields $\FF_q$ with $d(\FF_q)=-3$.
From the classification of unimodular hermitian modules we have the following table of orders of automorphism group
of unimodular irreducible hermitian modules over an imaginary quadratic extension $K$ of $\QQ$, with discriminant ${d}(K)=-3$.
\renewcommand{\arraystretch}{1.25}
\begin{center}\begin{tabular}{|c|l|}
\hline
${\rm dim}$ & $\#{\rm  Aut}$  \\
\hline
$2-5$& ---\\
$6$& $2^9\cdot 3^7 \cdot 5 \cdot 7$\\
$7$& ---\\[2pt]
$8$& $2^{14}\cdot 3^6 \cdot 5^2 \cdot 7$\\
$9$& $2^{8}\cdot 3^{12} \cdot 5 \cdot 7$\\
$10$& $2^{15}\cdot 3^9 \cdot 5^2 , \, 2^{17}\cdot 3^5 \cdot 5^2 \cdot 7$\\
\hline
\end{tabular}
\end{center}
\centerline{Table~1}
\begin{proposition}\label{d_3}
There is no optimal curve $C$ of genus $g(C)$ with $2 \le g(C) \le 10$ over a finite field $\FF_{q}$ of discriminant  $-3$ if $q \not = 3$.
\end{proposition}
\begin{proof}
We consider the case of an optimal curve $C$ over $\FF_{q}$.
From  Table~1 it follows that it cannot have genus $2$, $3$, $4$, $5$ and $7$,
since there is no corresponding irreducible unimodular hermitian module and hence no irreducible polarization.

If ${\rm  char}(\FF_{q})=p \ge 13$ and $g(C)=6,\;8,\;9$ or $10$ we would get that 
$\#{\rm  Aut}_{\FF_q}({\rm  Jac}(C))>2\cdot 84 \cdot (g-1)$ and hence $\#{\rm  Aut}{\FF_q}(C)>84 \cdot (g-1)$
which contradicts the upper bound mentioned above. 
\par
For characteristics $<13$ we can apply the Singh bound (cf.\ \cite{Singh}), viz. 
$$
\#{\rm  Aut}_{\bar{\FF_q}}(C) \le \frac{4pg^2}{p-1}\left(\frac{2g}{p-1}+1 \right)\left(\frac{4pg^2}{(p-1)^2}+1 \right).
$$
For the five characteristics which we have to examine, namely $2$, $3$, $5$, $7$ and $11$, the Singh bound is smaller than half the number
of automorphisms of the unimodular irreducible module, and hence the module cannot correspond to a polarization of a Jacobian.
\end{proof}
Combining these results with the non-existence of curves of defect $1$, i.\ e.\ with 
$\#C(\FF_q)=q+1\pm g[2\sqrt{q}] \mp 1$ (see \cite{Shabat}, page 22), we get the following theorem.
\begin{theorem}
If $C$ is a curve of genus $g$ with $3 \le g \le 10$ over a finite field $\FF_{q}$ of discriminant $-3$ and
$q \neq 3$ then 
$$
|\#C(\FF_{q})-q-1| \le g[2\sqrt{q}] - 2.
$$
\end{theorem}
\end{subsection}
\begin{subsection}{Optimal curves over finite fields with discriminant $-4$.}
In this subsection we examine the case of optimal curves of genus  $g \le 10$ over a finite field of the discriminant $d=-4$.
From the classification of unimodular hermitian forms with the discriminant $-4$ we have the following table 
of orders of automorphism groups of  irreducible hermitian modules
over an imaginary quadratic extension $K$ of $\QQ$ with discriminant $d(K)=-4$.
\renewcommand{\arraystretch}{1.25}
\begin{center}
\begin{tabular}{|c|l|}
\hline
${\rm dim}$&  $\#{\rm  Aut}$ \\
\hline
$2-3$& --- \\
$4$& $2^{10}\cdot 3^2 \cdot 5$ \\
$5$& --- \\[2pt]
$6$& $2^{15}\cdot 3^2 \cdot 5$ \\
$7$& $2^{11}\cdot 3^4 \cdot 5 \cdot 7$ \\
$8$& $2^{15}\cdot 3^5 \cdot 5^2 \cdot 7, 2^{22}
\cdot 3^2 \cdot 5 \cdot 7, 2^{21}\cdot 3^4 \cdot 5^2, 2^{15} \cdot 3^{2} \cdot 5 \cdot 7, 2^{21}\cdot 3^2$ \\
$9$& $ 2^{16}\cdot 3^{3} \cdot 5, 2^{10}\cdot 3^{4} \cdot 5^2 \cdot 7, 
2^{10}\cdot 3^{4} \cdot 5^2 \cdot 7$, \\
$10$& $ 2^{27}\cdot 3^{4} \cdot 5^2 \cdot 7, 2^{25}\cdot 3^{4} \cdot 5^2,
2^{17}\cdot 3^{5} \cdot 5 \cdot 7,2^{25} \cdot 3^2,$ \\
& $2^{19} \cdot 3^2\cdot 5 \cdot 7, 2^{17}\cdot 3^3 \cdot 5, 2^{11}\cdot 3^4 \cdot 5^2 \cdot7, 2^{23}\cdot 3\cdot 5,2^{19}\cdot 3^2,$\\
 & $2^{11}\cdot 3^4 \cdot 5^2, 2^{10}\cdot 3^4 \cdot 5^2, 2^{25} \cdot 3^3 \cdot 5$ \\[2pt]
\hline
\end{tabular}
\end{center}
\centerline{Table 2} 

\begin{proposition}
There is no optimal curve of genus $g$ with $2\le g \le 10$ over a finite field $\FF_{q}$ for
$d(\FF_{q})=-4$ and $q \neq 2$. 
\end{proposition}
\begin{proof}
Table~2 shows that
there is no such optimal curve of genus $2$, $3$ and $5$ since 
there are no corresponding irreducible unimodular hermitian modules of dimensions $2$, $3$ and $5$, respectively. 
\par
Let $C$ be an optimal curve of genus $g \le 10$ over a finite field $\FF_q$ of  discriminant $-4$. 
In case $ p \ge 13$ the unimodular hermitian module $(\cO_{K}^g,h)$ corresponding to
an optimal curve $C$ of genus $2 \le g \le 10$ should have $\#{\rm Aut}(\cO_{K}^g,h) \le 168\cdot(g-1)$,
but for $g \in \{4,6,7,8,9,10 \}$ there are no such hermitian modules as Table~2 shows.

It is easy to check that there are no $q=p^n$ with $p \in \{2,\; 3,\; 7,\; 11\}$ and with $d(\FF_q)=-4$.
For example, if $p=2$ and $n>1$ then $(m/2)^2= {2^n-1}\equiv {3}\; (\mbox{mod} \; 4)$,
 but $3$ is not a square ${\rm mod} \; 4$.
 
If ${\rm  char}(\FF_{q})=5$ one sees that number of the automorphisms of
the unimodular irreducible hermitian module is greater  than twice the Singh bound.
\end{proof}

Since there are no curves with defect $1$ (cf.\ \cite{Shabat}), we get the following theorem.
\begin{theorem}
If $C$ is a curve of genus $g$ with $3 \le g \le 10$ over a finite field $\FF_{q}$ with discriminant $-4$ and $q \neq 2$ then 
$$
|\#C(\FF_{q})-q-1| \le g[2\sqrt{q}] - 2.
$$
\end{theorem}
\end{subsection}
\begin{subsection}{Optimal curves over finite fields with discriminant $-7$.}
As before we begin with producing the table of orders of automorphism groups of the irreducible unimodular
hermitian modules over an imaginary quadratic extension $K$ of $\QQ$ with discriminant ${d}(K)=-7$.
\renewcommand{\arraystretch}{1.25}
\begin{center}
\begin{tabular}{|c|l|}
\hline
${\rm dim}$& $\#{\rm  Aut}$ \\
\hline
$2$& --- \\
$3$& $2^4\cdot3\cdot7$ \\
$4$& $2^{7}\cdot 3^2$ \\
$5$& $2^{8}\cdot 3 \cdot 5$ \\
$6$& $2^{9}\cdot 3^2 \cdot 7^2,\; 2^{10}\cdot 3^2 \cdot 5,\; 2^{9}\cdot 3^2
 \cdot 5,\; 2^{5}\cdot 3^2 \cdot 5 \cdot 7 $ \\
$7$& $ 2^{10}\cdot 3^4 \cdot 5 \cdot 7,\; 2^{11}\cdot 3^2 \cdot 5 \cdot 7,\; 2^{11}\cdot 3^3 \cdot 7,\;
 2^{8}\cdot 3^2 \cdot 5 \cdot 7,$\\
&$ 2^{11}\cdot 3^2,\; 2^{10}\cdot 3^2,\; 2^{6}\cdot 3^3 
\cdot 5,\; 2^{10}\cdot 3,\; 2^{7}\cdot 3 \cdot 5$ \\
$8$& $ 2^{14}\cdot 3^5 \cdot 5^{2} \cdot 7,\; 
2^{15}\cdot 3^2 \cdot 5 \cdot 7,\; 
2^{15}\cdot 3^4,\;
2^{12}\cdot 3^2 \cdot 5 \cdot 7,$\\

&$2^{9}\cdot 3^5 \cdot 5,\; 
2^{12}\cdot 3^2 \cdot 5,\; 
2^{14}\cdot 3^2,\;
2^{11}\cdot 3^2 \cdot 5,$\\

&$2^{11} \cdot 3^2,\; 
2^{9}\cdot 3^2 \cdot 5 \cdot 7,\; 
2^{12}\cdot 3^2 
2^{6}\cdot 3^2 \cdot 5 \cdot 7,$\\

&$2^{12} \cdot 3^2,\; 
2^{11}\cdot 3^2,\; 
2^{6}\cdot 3^2 \cdot 5 \cdot 7,\;,  
2^{10}\cdot 3^2,$\\

&$2^{7} \cdot 3^2 \cdot 5^{2},\; 
2^{13}\cdot 3,\; 
2^{8}\cdot 3^2 \cdot 5,\;,  
2^{8}\cdot 3^2,$\\

&$2^{9} \cdot 3^2,\;
2^{11},\; 
2^{8}\cdot 3^{4},\; 
2^{8}\cdot 3^4,\;,  
2^{7}\cdot 3^2,$\\

&$2^{9} \cdot 3,\;
2^{11},\; 
2^{7}\cdot 3^{2},\; 
2^{6}\cdot 3^2$ \\
\hline
\end{tabular}
\end{center}
\centerline{Table~3} 

Assume that there is an optimal curve $C$ of genus $g$ with 
$3 \le g \le 7$ over a finite field $\FF_{q}$. 
By the classification list and Proposition~\ref{proj} we see that there is a projection
$f_k:C \rightarrow E$ of degree $2$ on an optimal elliptic curve $E$, hence there exists an involution
$\sigma \in {\rm  Aut}_{\FF_{q}}$ so that $C/\langle \sigma \rangle \cong E$.
Let $G$ be a Sylow $2$-subgroup of the automorphism group ${\rm  Aut}_{\FF_{q}}(C)$ containing $\sigma$.
The involution
$\sigma$ does not belong to the center of the group $G$, since otherwise 
${\rm  Aut}_{\FF_{q}}(C/\langle \sigma \rangle) \cong {\rm  Aut}(E)=\{-1,1\}$ 
contains the factor group $G/\langle \sigma \rangle$ of order greater than $2^3$ as the table shows.
Since the center of $G$ is not trivial we can choose an involution $\tau$ in
the center of $G$. Then group generated by $\sigma$ and $\tau$
is isomorphic to $\ZZ/2\ZZ\times \ZZ/2\ZZ$ and 
hence we have the following diagram of coverings
$$
\xymatrix{
 &\ar[ld]_{2:1}C\ar[d]^{2:1}\ar[rd]^{2:1}&\\
 C/\langle \sigma \rangle \cong E\ar[rd]^{2:1}& C/\langle \sigma  \tau \rangle \ar[d]^{2:1} & C/\langle \tau \rangle \ar[ld]_{2:1} \\
&C/\langle \sigma, \, \tau \rangle &\\}
$$

Relations between idempotents in ${\rm End}({\rm Jac}(C))$ imply that
we have an isogeny
\begin{equation}\label{splitJac}
{\rm  Jac}(C)\times \left({\rm  Jac}(C/\langle \sigma, \, \tau \rangle)\right)^2 \sim E
\times {\rm  Jac}(C/\langle \sigma  \tau \rangle)\times {\rm  Jac}(C/\langle \tau \rangle),
\end{equation}
(cf.\ for example Thm.\ A of \cite{KR}).

\begin{lemma}\label{d-7g3hyperellitic}
An optimal curve of genus $3$ over the finite field $\FF_{q}$ is not a hyperelliptic.
\end{lemma}
\begin{proof}
Since there exists an irreducible unimodular hermitian module
of dimension $3$, there exists an irreducible principal polarization on the abelian variety $E^3$, where $E$
is  an optimal elliptic curve. Hence by \cite{Oort} and the Decent theory (or more precisely,
by the Theorem 4 and 5, in the Appendix in \cite{Lauter2}) 
there exists an optimal curve of genus $3$ over the finite
field $\FF_{q}$.
\par
If $C$ is a hyperelliptic optimal curve of genus $3$ over $\FF_q$ with  hyperelliptic involution $\iota$ (actually,
then $\iota=\tau$) then we have
$$
{\rm  Jac}(C)\sim E \times {\rm  Jac}(C/{\langle \sigma \tau \rangle}),
$$
where $\sigma$ an involution as constructed before.
Now the genus of the quotient curve $C/{\langle \sigma  \tau \rangle}$ is $2$, but there is no
optimal curve of genus $2$, and hence there is no hyperelliptic optimal curve of genus $3$.
\end{proof}

The key result of this subsection is the following proposition.
\begin{proposition}\label{optimald7}
There is no optimal curve $C$ of genus $g$ with $2 \le g \le 10$ and $g \not =3$ over a finite field $\FF_{q}$ of  discriminant $-7$.
\end{proposition}
\begin{proof}
If $C$ is an optimal curve of genus $4$, 
then by  Hurwitz's formula and by the nonexistence of optimal curves of genus $2$, the quotient curves $C/{\langle \tau \rangle}$ and $C/{\langle \sigma  \tau \rangle}$ 
have genus $<2$. Therefore,  the dimension of the left-hand-side of the equation~(\ref{splitJac}) 
is greater or equal to $4$, while the dimension of the right-hand-side is less or equal to $3$, 
which is impossible. 
\par
If $C$ is an optimal curve of genus $5$ then according to the choice of $\tau$ it follows
that $G/\langle \tau \rangle \subseteq {\rm  Aut}_{\FF_{q}}(C/\langle \tau \rangle)$ and
$\#G/\langle \tau \rangle \ge 2^6$. By Table~3 it follows that
$g(C/\langle \tau \rangle) \ge 4$ or $g(C/\langle \tau \rangle)=0$ and by Hurwitz's formula 
$g(C/\langle \tau \rangle)\le 3$. Hence the quotient curve $C/\langle \tau \rangle$ is a rational 
curve and $C/{\langle \sigma, \; \tau \rangle}$ is also rational.
Therefore the dimension of the left-hand-side of the equation~(\ref{splitJac}) is $5$ and the dimension
of the right-hand-side is less or equal to $2$, since an optimal curve of genus $3$ is not hyperelliptic.
This is a contradiction.
\par
If $C$ is an optimal curve of genus $6$ 
then  the dimension of the left-hand-side of the equation~(\ref{splitJac}) is an even number $\ge 6$,
but the dimension of the right-hand-side  $\le 7$
and cannot be equal to $6$, since there is no an optimal curve of genus $2$.
\par
If $C$ is an optimal curve of genus $7$, then
the quotient curve $C/ \langle \tau \rangle$ 
is either an optimal curve of genus $3$ or a rational curve, since $\# G/ \langle \tau \rangle \ge 2^4$. 
In the first case there is only one choice for the order of a Sylow $2$-subgroup G, namely $\#G=2^5$ as Table~3 shows. Then  $\#G/ \langle \tau \rangle =2^4$ and since $\#G/ \langle \tau \rangle \subset {\rm  Aut}_{\FF_{q}}(C)$ 
we get that the quotient curve $C/ \langle \tau \rangle$ is hyperelliptic by a similar argument, but there is no  hyperelliptic optimal curve of genus $3$ by Lemma~\ref{d-7g3hyperellitic}. 
In the second case,  $\tau$ is the hyperelliptic involution and $\langle \tau, \, \sigma \rangle \cong \ZZ/2\ZZ \times \ZZ/2\ZZ$. Now (\ref{splitJac}) gives the following splitting of Jacobian (up to isogeny)
$$
{\rm  Jac}(C)\sim
{\rm  Jac}(C/{\langle \sigma \rangle}) \times {\rm  Jac}(C/{\langle \sigma  \tau \rangle}).
$$
We know that ${\rm  Jac}(C/{\langle \sigma \rangle})$ has the dimension $1$.
Therefore $C/{\langle \sigma  \tau \rangle}$ has to be an optimal curve of genus $6$, but such a curve does  not exist as we just showed.

By exactly the same arguments it follows that there are no optimal curves of genus $8.$
\end{proof}

As a consequence of Proposition~\ref{optimald7} we obtain the following theorem.
\begin{theorem}
If $C$ is a curve of genus $g$ with $4 \le g \le  8$ over a finite field $\FF_{q}$ of discriminant $-7$ then 
$$
|\#C(\FF_{q})-q-1| \le g[2\sqrt{q}] - 2.
$$
\end{theorem}
\end{subsection}
\begin{subsection}{Optimal curves over the finite field $\FF_{q}$ with discriminant $-8$.}
We start by giving the table of orders of automorphism groups of the irreducible hermitian modules 
over an imaginary quadratic extension $K$ of $\QQ$ with discriminant $d(K)=-8$.

\renewcommand{\arraystretch}{1.25}
\begin{center}
\begin{tabular}{|c|l|}
\hline
${ \rm dim}$&  $\#{\rm  Aut}$ \\
\hline
$2$& $2^4\cdot 3$ \\
$3$& --- \\
$4$& $2^{7}\cdot 3,\; 2^{8}\cdot 3^2, 2^{9}\cdot 3^2$ \\
$5$& $2^{5}\cdot 3^2\cdot5 $ \\
$6$& $2^{10}\cdot3^2\cdot5,\; 2^{12}\cdot3^3, 2^8\cdot3^4\cdot5,\; 2^{13}\cdot3^4$\\
&$2^{10}\cdot3^2 \cdot 5,\; 2^{11}\cdot3^2,\; 2^{10}\cdot3,\; 2^6\cdot 3^2\cdot 5,\; 
2^7\cdot 3^2,\; 2^{11},\; 2^{12}\cdot3$ \\
$7$& $2^{10}\cdot3^4\cdot5\cdot7,\; 2^{10}\cdot3^2\cdot5,\; 
2^9\cdot3^3\cdot5,\; 2^{5}\cdot3^2 \cdot 5 \cdot 7,$\\
&$ 2^{8}\cdot3^2,\; 2^{10}\cdot3^2,\; 2^{10}\cdot3,\;
2^5\cdot 3^2\cdot 5,\; 2^8\cdot 3,\; 2^{5}\cdot 3^3,\; 2^{5}\cdot3^4$ \\
\hline
\end{tabular}
\end{center}
\centerline{Table 4} 

\begin{proposition}
There is no optimal curve of genus $g$
with $3\le g \le 7$ over a finite field $\FF_q$ of discriminant $-8$ and ${\rm char}(\FF_q) \not =3$ . 
\end{proposition}
\begin{proof}
There is no optimal curve of genus $3$ over a finite field $\FF_q$ of discriminant $-8$, 
since does not exist a corresponding hermitian module.
\par
Assume that $C$ is an optimal curve  of genus $4$ over $\FF_{q}$ of  discriminant 
$-8$ and ${\rm  char}(\FF_{q})>5$.
If  $C$ is hyperelliptic then by Table~4 we see that  $\#{\rm  Aut}_{\FF_{q}}(C)=\#{\rm  Aut}_{\FF_{q}}({\rm  Jac}(C)) > 84\cdot (4-1)$, a contradiction to the bound we observed. If $C$ is not hyperelliptic then Tabel~4 yields 
that the automorphism group ${\rm  Aut}_{\FF_{q}}(C)$ has a subgroup $G$ of order $2^6$. 
The group $G$ has a normal subgroup $H$ of order $2$ and the automorphism group of the quotient 
curve $C/H$ has a subgroup isomorphic to $G/H$, but this is impossible 
if $g(C)>0$, since this order contradicts Table~4.
Therefore $C/H$ is a rational curve and $C$ is hyperelliptic.
\par
If $C$ is an optimal curve of genus $5$ and ${\rm  char}(\FF_{q}) \ge 7$,
then by Table~4 $\#{\rm  Aut}_{\FF_{q}}(C) \ge \frac{1}{2}\#{\rm  Aut}_{\FF_{q}}({\rm  Jac}(C)) >84\cdot(5-1)$, 
and this contradicts the automorphism bound.
\par
The case of genus $6$  is similar to that of genus $5$.
\par
Now we consider optimal curves of genus $7$.
Assume that $C$ is an optimal curve of genus $7$ over $\FF_{q}$ of discriminant $-8$.
We split the set of orders of the automorphism groups of the irreducible unimodular hermitian modules of 
the dimension $7$ into three subsets, $T_1=\{2^{10}\cdot3^4\cdot5\cdot7,\, 
2^{10}\cdot3^2\cdot5,\, 2^9\cdot3^3\cdot5,\,
2^9\cdot3^3\cdot5,\, 2^5\cdot3^2\cdot5\cdot7,\, 2^5\cdot3^4, \, 2^5\cdot3^2\cdot5\}$,  
$T_2=\{2^8\cdot3^{2},\, 2^{10}\cdot3^2,\, 2^{10}\cdot3,\, 2^8\cdot3 \}$  and 
$T_3=\{2^5\cdot3^3\}$.

If ${\rm  char}(\FF_{q})\ge 11$ and $\#{\rm  Aut}_{\FF_q}({\rm  Jac}(C)) \in T_1$ then we get that
$\#{\rm  Aut}_{\FF_q}(C) \ge \frac{1}{2}\,\#{\rm  Aut}_{\FF_q}({\rm  Jac}(C))>84(7-1)$, and 
hence there does not exist an optimal curve in this case.
\par
If $\#{\rm  Aut}_{\FF_q}({\rm  Jac}(C)) \in T_2$ then
it is not hyperelliptic, since otherwise one can get a contradiction as before. The automorphism group
${\rm  Aut}_{\FF_q}({\rm  Jac}(C))$ contains a subgroup $G$ of order~$2^7$, 
which has a normal subgroup $H$ of order~$2$. The automorphism group of the quotient 
curve $C/H$ then has a subgroup isomorphic to $G/H$. 
Since there is no optimal curve of genus $3$, $4$, $5$, $6$ and~$2^6$ does not divide 
the orders of automorphism groups of optimal curves of genus $1$ and $2$, it follows that
$C/H$ is a rational curve. But this is also impossible, since  $C$ is not hyperelliptic. 
\par
Let $\#{\rm  Aut}_{\FF_q}({\rm  Jac}(C)) \in T_3$ and $G$ be a Sylow $2$-subgroup of
${\rm  Aut}_{\FF_{q}}(C)$. 
Then similarly as for discriminant $-7$ there are  two 
involutions $\sigma,\, \tau \in G$ with the same properties.
Here we would like to remark that there are two cases: the first hermitian form 
give us an involution $\sigma$ by Proposition~\ref{proj}. In the second case, this does not work, but if we apply
an isomorphism given in terms of standard bases as $e_i \rightarrow e_i$ if $i \not = 4$ and 
$e_4 \rightarrow e_4-e_5+e_6$, then by applying Proposition~\ref{proj} we again have an involution $\sigma$.
Therefore there is an isogeny
$$
{\rm  Jac}(C)\times {\rm  Jac}(C/\langle\sigma,\, \tau \rangle)^2 \sim {\rm  Jac}(C/\langle \sigma \rangle)
\times {\rm  Jac}(C/\langle \sigma  \tau \rangle)\times {\rm  Jac}(C/\langle \tau \rangle).
$$
The dimension of the left-hand-side is $>6$, whereas the dimension of the right-hand-side is $<6$,
because one of the factors is an elliptic curve and the other two have genus~$\le 2$.
\par
Finally, note that there is no finite field $\FF_q$ of ${\rm char}(\FF_q)\in \{2,5,7\}$
of discriminant~$-8$. 
Indeed, $2(4^n-1)\equiv 2\;(\mbox{mod}\; 4)$ is not square in $({ \rm mod}\; 4)$,
similarly, $5^{2n+1}-2$ is not a square $({ \rm mod}\; 5)$ and 
$7^{2n+1}-2$ is not a square $({ \rm mod}\; 3)$.
\end{proof}

\begin{theorem}
If $C$ is a curve of genus $g$ with $3 \le g  \le 7$ over the finite field $\FF_{q}$ of discriminant $-8$ and
${\rm  char}(\FF_{q}) \neq 3$  then 
$$
|\#C(\FF_{q})-q-1| \le g[2\sqrt{q}] - 2.
$$
\end{theorem}
\end{subsection}
\end{section}
\begin{section}{Optimal curves over finite fields with discriminant $-11$.}
\begin{subsection}{Optimal elliptic curves and optimal curves of genus $2$.}
In this subsection we study optimal elliptic curves over $\FF_{q}$, their
$\bar{\FF}_{q}$--isomorphism classes and the $\FF_{q}$--isomorphism classes in an isogeny class.
In addition, we give the results of concrete calculations for finite fields $\FF_{q}$ with discriminant $-11$ and $ q <  10^4$.
As we mentioned in Section~1, with that information we can find the number of isomorphism classes
of abelian varieties which are candidates for being Jacobians of optimal curves. 
\begin{proposition}\label{elliptic11}
Up to isomorphism over a finite field $\FF_{q}$ with discriminant $-11$
there exists exactly one maximal (resp.\ minimal) optimal elliptic curve $E$ over $\FF_{q}$.
\end{proposition}
\begin{proof}
Let $E$ be an optimal elliptic curve over $\FF_q$.
The trace of Frobenius of the optimal elliptic curve
over ${\FF}_{q}$ equals $m=\pm [2\sqrt{q}]$ and its 
Frobenius endomorphism $F$ satisfies the equation $F^2\pm m F + q=0$.
This fixes the isogeny class of the optimal maximal (resp.\ minimal)
elliptic curve. Then ${\rm   End}_{{\FF}_{q}}(E)$ contains a ring
${\ZZ}[X]/(X^2\pm mX+q)$, which is the ring of integers $\cO_K$
of $K={\QQ}(\sqrt{-11})$, hence ${\rm   End}_{{\FF}_{q}}(E)=\cO_K$.
Since the class number of $\cO_K$ is $1$ the ${\FF}_{q}$-isogeny 
class of $E$ contains one $\bar{\FF}_{q}$-isomorphism class
and such an $\bar{\FF}_{q}$-isomorphism class contains exactly two ${\FF}_{q}$-isomorphism classes, since 
${\rm   Aut}_{\bar{\FF}_{q}}(E)={\ZZ}/2{\ZZ}$ (see the "mass formula" in \cite{Geer}). 
One of these classes corresponds to a maximal optimal curve, while its twist then gives a minimal elliptic curve. 
\end{proof}
\begin{example}
In the following table we give for various $\FF_q$ of discriminant $-11$ an optimal
maximal and an optimal minimal elliptic curve. In the table 
an entry $[a,b]$ means the elliptic curve over $\FF_q$ given by 
an equation $y^2=x^3+a\,x+b$.
\renewcommand{\arraystretch}{1.25}
\begin{center}
\begin{tabular}{|c|l|l|}
\hline
${q}$& maximal& minimal\\
\hline
${23}$& $[1,11]$& $[12,8]$\\
${59}$& $[2,22]$& $[30,47]$\\
${113}$& $[5,44]$& $[101,10]$\\
${383}$& $[1,91]$& $[46,301]$\\
${509}$& $[1,70]$& $[382,136]$\\
${563}$& $[2,363]$& $[282,538]$\\
${1193}$& $[5,11]$& $[1061,619]$\\
${1409}$& $[6,221]$& $[940,1365]$\\
${3083}$& $[2,1009]$& $[1542,2053]$\\
${4973}$& $[1,1748]$& $[3730,2705]$\\
${6323}$& $[2,221]$& $[3161,818]$\\
\hline
\end{tabular}
\end{center}
\end{example}
\begin{remark}
The points of order $2$ of the elliptic curve consist of two
Galois orbits, one of length $1$ and one of length $3$. Indeed,
In view of $-11=m^2-4q$ and $\#E(\FF_q)=q+1\pm m$ we see that
$m$ is odd, hence $\#E(\FF_q)$. This implies that there $\#E[2](\FF_q)=1$.
\end{remark}
The case of optimal curves of genus $2$ was dealt with by J-P.\ Serre \cite{Shabat}. 
One finds that over the field $\FF_{q}$ with discriminant $-11$ an optimal curve of genus $2$ always exists. 
J.-P.\ Serre constructed such curves by gluing two optimal elliptic curves. 
Here we show how an optimal curve of genus $2$ over $\FF_q$ can be obtained 
by gluing two optimal elliptic curves and that there is only one $\FF_{q}$-isomorphism class of such curves. 
It also follows that up to $\FF_q$-isomorphism there are at most two degree $2$ maps 
of an optimal curve of genus $2$ to a fixed optimal elliptic curve.
\begin{proposition}\label{genustwo11}
Up to isomorphism over the field ${\FF}_{q}$ of discriminant $-11$
there exists exactly one maximal  (resp.\ minimal) optimal curve
$C$ of genus $2$ over ${\FF}_{q}$, viz., the fibered product over ${\PP}^1$
of the two maximal (resp.\ minimal) optimal elliptic curves
$$
y^2=f(x) \quad  \mbox{and} \quad  y^2=f(x)(\alpha\, x+\beta),
$$
where $f(x)$ is polynomial of degree $3$.
\end{proposition}
\begin{proof} 
By Serre we know that exists at least one maximal (minimal) optimal curve of genus $2$ over
 $\FF_q$with discriminant $-11$.
Let $C$ be an optimal curve of genus $2$ over ${\FF}_{q}$.
Using the equivalence of categories given in Section~2 the theta divisor
of the canonical lift Jacobian ${\rm   Jac}(C)$ corresponds to a unimodular hermitian module
that has the unimodular irreducible hermitian form given by
$$
\left(
\begin{array}{cc}
2 & \frac{-1+\sqrt{-11}}{2}\\
\frac{-1-\sqrt{-11}}{2} & 2\\
\end{array}
\right).
$$
By Proposition~\ref{proj} it follows that $C$ admits a map $\psi$ of degree $2$
onto an optimal elliptic curve $E$. This means that the automorphism group
of $C$ contains two commuting automorphisms, one $\alpha$ 
corresponding to $\psi$ and the other one the hyperelliptic involution $\iota$.  It follows that $C$ is 
(the normalization of) a fibered product over 
$\PP^1=C/\langle \alpha,\iota\rangle$ of the two optimal elliptic curves
$C/\langle \alpha\rangle$ and $C/\langle\iota\alpha \rangle$. 
Now we treat the maximal case; the minimal case follows then by twisting.

Let $C \cong E \times_{\PP^1}E^{\prime}$ be such a fibered product
of two maximal elliptic curves. Of course $E^{\prime}$ is
isomorphic to $E$, the curve given in \ref{elliptic11}.
In order that the fibered product of $E$ and $E^{\prime}$ over the $x$-line 
${\PP}^1$ is of genus $2$, three of their ramification points must coincide.
If $E$ is given by an equation of the form $y^2=f(x)$ then the polynomial $f(x)$ is irreducible over ${\FF}_{q}$
and it follows that $E^{\prime}$ can be written as $y^2= (\alpha\, x +\beta)\, f(x)$
for $\alpha, \beta \in {\FF}_{q}$ with $\alpha \neq 0$. By the result of Proposition~\ref{elliptic11}
there exists an automorphism $g \in {\rm   PGL}(2,{\FF}_{q})$ of
${\PP}^1$ that permutes the roots of $f(x)$ and sends $\infty$ to $x=-\beta/\alpha$.
Since the roots of $f(x)$ generate a cubic extension of 
${\FF}_{q}$ the element $g$ is necessarily of order $3$ or $1$, hence there are at
most two such equations for the curve $E^{\prime}$ that can occur. 
\end{proof}
\begin{remark}\label{genus2_degree6}
From the fact that an optimal curve $H$ of genus $2$ over $\FF_q$ with $d(\FF_q)=-11$ 
is the fibered product of two optimal elliptic curves
it follows that $H$ can be given by an equation $z^2=F(x)$, where $F(x)$ is a polynomial of degree $6$.
\end{remark}
\begin{remark}\label{weierstrass11}
Since the points of order $2$ of the elliptic curve consist of two
Galois orbits, one of length $1$ and one of length $3$, it follows that
the ramification points of our  curve of genus~$2$ form two Galois orbits
of length~$3$.
\end{remark}
\begin{example}
Here we give a list of maximal optimal curves of genus $2$ as given Proposition~\ref{genustwo11}.
In the table we denote a maximal curve of genus $2$ given by an equation $z^2=\alpha\, x^6+\beta\, x^4+\gamma\,x^2+\delta$  by $(\alpha,\beta,\gamma,\delta)$ 
and a maximal elliptic curve given by an equation $y^2=(x^3+a\,x+b)(c\, x+d)$ by $[a,b;c,d]$
Minimal optimal curves of genus $2$ are obtained by twisting.
\renewcommand{\arraystretch}{1.25}
\begin{center}
\begin{tabular}{|c|l|l|}
\hline
$q$& $[a,b;c,d]$& $(\alpha,\beta,\gamma,\delta)$\\
\hline
${23}$& $[1,11;1,19]$& $(1,12,3,10)$\\[2pt]
${59}$& $[2,22;1,49]$& $(1,30,7,39)$\\[2pt]
${113}$& $[5,44;1,24]$& $(1,41,38,112)$\\[2pt]
${383}$& $[1,91;1,135]$& $(1,361,290,356)$\\[2pt]
${509}$& $[1,70;2,208]$& $(191,431,191,501)$\\[2pt]
${563}$& $[2,363;1,189]$& $(1,559,195,212)$\\[2pt]
${1193}$& $[5,11;1,1017]$& $(528,1074,1072,657)$\\[2pt]
${1409}$& $[6,221;3,1135]$& $(835,187,516,278)$\\[2pt]
${3083}$& $[2,1009;1,2569]$& $(1,1542,259,1880)$\\[2pt]
${4973}$& $[1,1748;2,2341]$& $(1865,987,4365,4633)$\\[2pt]
${6323}$& $[2,221;1,1274]$& $(1,2501,520,3218)$\\[2pt]
\hline
\end{tabular}
\end{center}
\end{example}
\end{subsection}
\begin{subsection}{Optimal curves of genus $3$, $4$ and $5$.}
As J.-P.\ Serre showed (cf.\ \cite{Serre}) there is no optimal curve of genus~$3$ over a field $\FF_{q}$ with
discriminant~$-11$, since there is no irreducible principal polarization of the Jacobian. 
Therefore we may start with genus~$4$.
\begin{lemma}
If there is an optimal curve $C$ of genus $4$ over  $\FF_{q}$ with discriminant $-11$,
then there exist two commuting involutions $\sigma, \tau \in {\rm Aut}_{\FF_q}(C)$ such that
$C/\langle \sigma \rangle$ is an optimal curve of genus $2$ over $\FF_q$ and $C/\langle \tau \rangle$ and $C/\langle \sigma \tau \rangle$ are optimal elliptic curves over $\FF_q$.
Moreover, we can find irreducible polynomials $h_1(x), h_2(x) \in \FF_q[x]$ of degree $3$ such that
$C/\langle \sigma \rangle$ is  given by an equation $y^2=h_1(x)\cdot h_2(x)$ and $C/\langle \tau \rangle$, $C/\langle \sigma \tau \rangle$ are given by equations $y^2=h_1(x)$ and $y^2=h_2(x)$.
\end{lemma}
\begin{proof}
Assume that there is an optimal curve $C$ of genus $4$ over $\FF_q$ with discriminant $-11$.
From the classification of unimodular
hermitian modules it follows that $\# {\rm  Aut}_{\FF_{q}}({\rm  Jac}(C)) \in \{2^5\cdot 3^2, \, 2^4\cdot 3\cdot 5, \, 2^4\cdot3^2 \}$. 
Let $G$ denote a Sylow $2$-subgroup of  ${\rm  Aut}_{\FF_{q}}(C)$ and $\sigma \in G$ an involution in the center of $G$.
Since $G/\langle \sigma \rangle \subseteq {\rm  Aut}_{\FF_{q}}(C/\langle \sigma \rangle)$ it follows that
$C/\langle \sigma \rangle$ is either a rational curve or an optimal curve of genus~$2$
since the order $\#{\rm Aut}(C/\langle \sigma \rangle)$ 
is too large for $C/\langle \sigma \rangle$ being an elliptic curve. Moreover,  
Proposition~\ref{proj} together with the classification list show that there is another involution $\tau \in G$ such that the quotient curve $C/\langle \tau \rangle$ is an optimal elliptic curve.
Therefore $ \sigma, \, \tau$ generate a group isomorphic to $ \ZZ/2\ZZ\times \ZZ/2\ZZ$.
\par
If $\sigma$ is the hyperelliptic involution of $C$, i.e.\ $C/\langle \sigma \rangle \cong \PP^1$, then
$$
{\rm  Jac}(C)\sim {\rm  Jac}(C/\langle \tau \rangle)\times {\rm  Jac}(C/\langle \sigma\,\tau \rangle).
$$ 
However, the dimension of the left-hand-side is $4$ and of that the right-hand-side is less or equal to $3$, 
so $\sigma$ cannot be the hyperelliptic involution.
\par
Now assume that $C/\langle \sigma \rangle$ is an optimal curve of genus $2$. 
Then there is the isogeny of Jacobians
$$
{\rm  Jac}(C)\times {\rm  Jac}(C/\langle \sigma,\, \tau \rangle)\sim{\rm  Jac}(C/\langle \sigma \rangle)
\times {\rm  Jac}(C/\langle  \tau \rangle)\times {\rm  Jac}(C/\langle \sigma   \tau \rangle),
$$
again by comparing dimensions and since $g(C/\langle \sigma , \tau \rangle)=0$ (and there is no double
covering of a maximal elliptic curve by a maximal elliptic curve) we have that
$C/\langle \sigma  \tau \rangle$ is an optimal elliptic curve.
\par
We know by Remark~\ref{weierstrass11} that the ramification points of $C/\langle \sigma \rangle$ form two
Galois orbits, say $\{p_1, p_2, p_3 \}$ and $\{p_4, p_5, p_6 \}$.
We denote the inverse image of $p_i$ on $C$ by $\{ p_{i}^{\prime}, p_{i}^{\prime \prime}\}$. Then 
$\sigma( p_{i}^{\prime})=p_{i}^{\prime \prime}$ for all $i$, and we have $\tau(p_{i}^{\prime})=p_{i}^{\prime}$ 
and $\tau(p_{i+3}^{\prime})=p_{i}^{\prime \prime}$ for $i=1,2,3$. The involution $\sigma$ has two ramification
points, say $q_1$, $q_2$ and these are interchanged by $\tau$. We let the image of $q_1$ (and of $q_2$) on 
$C/\langle \sigma,\, \tau \rangle \cong \PP^1$ be equal to $\infty$ and we denote the image of $p_i$ on 
$C/\langle \sigma,\, \tau \rangle$ by $P_i$. Then it is clear that the branch points of 
$C/\langle  \tau \rangle \rightarrow C/\langle \sigma,\, \tau \rangle$ are $\{P_1, P_2, P_3, \infty \}$,
while those of $C/\langle  \sigma \tau \rangle \rightarrow C/\langle \sigma,\, \tau \rangle$ are
$\{P_4, P_5, P_6, \infty \}$. This shows that we can choose an equation $y^2=h_1(x)h_2(x)$ for 
$C/\langle  \sigma \rangle$ and equations for $C/\langle   \tau \rangle$ and $C/\langle  \sigma \tau \rangle$
as indicated in the statement.
\end{proof}
To check for the existence of an optimal curve of genus $4$ over $\FF_q$ with $d(\FF_q)=-11$,
it suffices by this lemma to check in view of Proposition~\ref{genustwo11}  whether the unique
optimal curve of genus $2$ over $\FF_q$ is of the form 
$y^2=h_1(x)h_2(x)$ with $h_1(x),\, h_2(x)$ as in the lemma.
This can easily be checked and we did this on a computer for $q <10^4$ with $p \not = 3$.
\begin{corollary}
There is no optimal curve of genus $4$ over a finite field $\FF_{q}$ of discriminant $-11$  if  $q <  10^4$ and 
${\rm char}(\FF_q) \not =3$.
\end{corollary}
Now we deal with genus $5$.
\begin{theorem}
There is no optimal curve $C$ of genus $5$ over a finite field $\FF_{q}$ of discriminant $-11$ 
and ${\rm  char}(\FF_{q})>5$.
\end{theorem}
\begin{proof}
Assume that $C$ is an optimal curve of genus $5$ over $\FF_{q}$ with discriminant $-11$.
From the classification list of the unimodular irreducible hermitian modules 
it follows that 
$\#{\rm  Aut}_{\FF_{q}}({\rm  Jac}(C))\in \{2^3\cdot3\cdot5\cdot11,\; 2^4\cdot 3\cdot 5,\; 2^4\cdot3\}$. 
\par
If ${\rm  char}(\FF)_{q}>5$ then $\#{\rm  Aut}_{\FF_{q}}({\rm  Jac}(C)) \not =2^3\cdot3\cdot5\cdot11$,
since $2^3\cdot3\cdot5\cdot11 \ge 84(5-1)$.
\par
If $\#{\rm  Aut}_{\FF_{q}}({\rm  Jac}(C)) \in \{2^4\cdot 3\cdot 5,\; 2^4\cdot3 \}$ 
then as for genus $4$  there are two involutions $\sigma$ and $\tau$  such that
$C/\langle \sigma \rangle$ is an optimal elliptic curve, 
$\sigma   \tau =\tau \sigma$ and $C/\langle \tau \rangle$ is not an optimal elliptic curve. 
As a result $C$ is the fibered product of
$C_1:=C/\langle \sigma \rangle$ and $C_2:=C/\langle \sigma   \tau \rangle$
over $C/\langle \sigma,\, \tau \rangle$ and hence we have the isogeny of Jacobians
$$
{\rm   Jac}(C/\langle \sigma   \tau \rangle)\times {\rm   Jac}(C/\langle \tau \rangle) \times E\sim
{\rm   Jac}(C/\langle \sigma, \tau \rangle)^2 \times{\rm   Jac}(C).
$$
In view of the non-existence of an optimal curve of genus $3$, the dimension of the left-hand-side 
$\leq 5$, but the dimension of the right-hand-side  $\geq 5$. The only possibility is that
the curves $C/\langle \sigma  \tau \rangle$, 
$C/\langle \tau \rangle$ have genus~$2$ and that $C/\langle \sigma, \tau \rangle$ has genus $0$.
But then  $C/\langle \sigma  \tau \rangle$ and $C/\langle \tau \rangle$ share two of their six ramification points and this contradicts \ Remark~\ref{weierstrass11}, 
which says that the six ramification points decompose in two Galois orbits of length $3$.
\end{proof}
\end{subsection}
\end{section}
\begin{section}{Optimal curves over finite fields with discriminant $-19$.}
\begin{subsection}{Optimal elliptic curves and optimal curves of genus $2$.}
In this subsection we explore optimal elliptic curves over $\FF_{q}$ and produce 
the results of calculations for finite fields $\FF_{q}$ with discriminant $-19$ and  $q <  10^3$.
 
Arguments similar to that for discriminant $-11$ give the following results.
\begin{proposition}\label{elliptic19}
Up to isomorphism over $\FF_{q}$ with discriminant $-19$ there exists exactly one maximal (resp. minimal) optimal elliptic curve $E$ over $\FF_{q}$.
\end{proposition}
\begin{remark}
The points of order $2$ of the elliptic curve consist of two
Galois orbits, one of length~$1$ and one of length~$3$.
\end{remark}
\begin{proposition}\label{genustwo19}
Up to isomorphism over the field ${\FF}_{q}$ with discriminant $-19$
there exists exactly one maximal  (resp.\ minimal) optimal curve
$C$ of genus $2$ over ${\FF}_{q}$, viz., the fibered product over ${\PP}^1$
of two maximal (resp.\ minimal) optimal elliptic curves
of the form
$$
y^2=f(x) \quad {\rm   and} \quad  y^2=f(x)(\alpha\, x+\beta).
$$
\end{proposition}
\begin{remark}\label{weierstrass2}
As before (for discriminant $-11$), it follows that
the ramification points of our genus $2$ curve form two Galois orbits
of length $3$.
\end{remark}
\end{subsection}
\begin{subsection}{Optimal Curves of Genus $3$ and $4$.}
\begin{proposition}
Up to isomorphism over $\FF_{q}$ with discriminant $-19$ there is an optimal curve $C$ of
genus $3$ over $\FF_{q}$, namely the double covering of a (maximal or minimal) optimal elliptic curve.
\end{proposition}
\begin{proof}
Applying the equivalence of categories as described in Section~2 along 
with the classification of unimodular hermitian modules, we see that 
the canonical lift of the Jacobian of an optimal curve of genus $3$
has a polarization that corresponds to the unimodular hermitian form 
$$
\left(
\begin{array}{ccc}
2 & 1 &-1 \\
1 & 3 & (-3+\sqrt{-19})/2\\
-1 & (-3-\sqrt{-19})/2&3\\
\end{array}
\right)_.
$$
By \cite{Oort} this is the Jacobian of a curve, hence there exists a (minimal or maximal) optimal curve $C$ of genus $3$ over a finite field $\FF_q$
and $f_{1} : C \to E$ is a double covering onto an optimal elliptic curve by Proposition~\ref{proj} . 

To finish the proof we remark that
there is only one class of such irreducible unimodular hermitian 
modules, hence by the equivalence of categories and Torelli's theorem 
the optimal curve of genus $3$ is unique up to isomorphism over $\FF_{q}$.
\end{proof}
In the following we use the explicit description of the
automorphism groups of the irreducible unimodular
hermitian modules \cite{Schulze-letter} to prove the following theorem:
\begin{theorem}
There is no optimal curve of genus $4$ over a finite field $\FF_{q}$ with discriminant
$-19$ if $q <  10^3$ and $q\equiv 1 \, (\bmod\, 5).$ 
\end{theorem}
\begin{proof}
Let $C$ be an optimal curve of genus $4$ over $\FF_{q}$. 
The classification of the unimodular hermitian modules of rank $4$ and
 discriminant $-19$ tell us that there are $9$ 
irreducible unimodular hermitian modules up to isomorphism.
If a unimodular hermitian module corresponds to the Jacobian of a curve $C$, 
then its automorphism group is isomorphic to the automorphism group of the 
polarized abelian variety ${\rm   Jac}(C)=(E^g, \Theta)$, hence by the
Theorem of Torelli the automorphism group
of the hermitian lattice and the automorphism group of the curve $C$
coincide if the curve is hyperelliptic; in the other case their
orders differ by a factor $2$.

From the structure of the automorphism group of the unimodular hermitian modules we get that an automorphism group 
${\rm  Aut}_{\FF_{q}}(C)$ of an optimal curve $C$ has subgroup isomorphic to either $\ZZ/2\ZZ\times \ZZ/2\ZZ$
or $\ZZ/5\ZZ$ (cf. the Appendix).

To manage the case  when ${\rm  Aut}(C)$ has subgroup $\ZZ/5\ZZ$ 
(this holds for the hermitian unimodular modules in the cases $(2)-(4)$) 
we  consider the cyclic Galois covering $C \to C/\left(\ZZ/5\ZZ \right)$, which is a Kummer cover, since  $q\equiv 1 \, (\bmod\, 5)$. 
Hurwitz's formula implies  that the quotient curve has genus $g'=0$ since $6=10\, (g' -1)+4\, t$ with $t$
the number of branch points, hence  $C/\left(\ZZ/5\ZZ \right)\cong \PP^1$. The number of rational points 
on the curve $C$ is $5s+l$, where $l$ is the number of rational ramification points. 
\par
If $d(\FF_q)=-19$ then
there are only three possibilities for the value $\#C(\FF_{q})$  $(\mbox{{\rm  mod}}\; 5)$, namely
 $\#C(\FF_{q})\equiv 0,\;1\; \mbox{or} \; 2 \; (\mbox{{\rm  mod}} \; 5)$.
\par
If $\#C(\FF_{q})\equiv 2\; (\mbox{{\rm  mod}}\; 5)$ then
two rational points must be ramification points since $0 \le l \le 4$. 
We may assume that the points $0$ and $\infty$ are ramified in the covering, 
and hence the curve $C$ can be given by an equation 
$$
z^5=\gamma x^{\nu_1}(x^2+ax+b)^{\nu_2},
$$
where $1 \le \nu_1, \nu_2 \le 4$ and the polynomial $x^2+ax+b$ is irreducible in $\FF_{q}[x]$. 
Since infinity is ramified, it follows that $\nu_1+2\nu_2 \neq 5$. 
If $\#C(\FF_{q})\equiv 1\; (\mbox{{\rm  mod}}\; 5)$ then
one rational point must be ramification point. 
One may assume that $\infty$ is ramified in the covering, and hence the curve $C$ can be given by an equation
$$
z^5=\gamma (x^3+ax+b)^{\nu},
$$
where $1 \le \nu \le 4$ and the polynomial $x^3+ax+b$ has to be irreducible in $\FF_{q}[x]$.
If $\#C(\FF_{q})\equiv 0 \; (\mbox{{\rm  mod}}\; 5)$ then
there are no rational ramification points. Therefore we may assume that
the curve $C$ is given by equation
$$
z^5=\gamma (x^2+ax+b)^{\nu_1}(x^2+c)^{\nu_1},
$$
where the polynomials $x^2+ax+b$ and $x^2+c$ are irreducible over $\FF_{q}$.
Since infinity is unramified, it follows that $2\nu_1+2\nu_2 \equiv\; 0\; (\mbox{{\rm  mod}}\; 5)$.  
A computer search showed that there are no optimal curves of these forms. This finishes the cases when
${\rm Aut}(C)$ contains $\ZZ/5\ZZ$.
\par
For the other cases in view of the Theorem of Torelli  we treat the hyperelliptic and 
the non-hyperelliptic case separately.
\par
If $C$ is a hyperelliptic curve with the hyperelliptic involution $\iota$, then by Theorem of Torelli the automorphism group ${\rm  Aut}(C)$ equals to ${\rm  Aut}({\rm  Jac}(C))$. 
The automorphism groups of the unimodular hermitian module 
contain an involution $\neq \iota$ (cf. the Appendix), and hence there is a non-hyperelliptic involution $\sigma $
in ${\rm  Aut}(C)$. Therefore we can consider the following diagram:
$$
\xymatrix{
&\ar[ld]_{2:1}C\ar[rd]^{2:1}&\\
 H_1:=C/\langle\sigma\rangle\ar[rd]^{2:1}& & C/\langle\iota\sigma\rangle=:H_2 \ar[ld]_{2:1} \\
&C/\langle\sigma, \iota\sigma\rangle&\\}
$$
The quotient curve $C/\langle\sigma,\iota\sigma \rangle$ is a rational line, 
since it is a quotient of $C/\iota \cong\PP^1$. 
In view of the isogeny ${\rm  Jac}(C) \sim {\rm  Jac}(C/\langle \iota\sigma\rangle) \times {\rm  Jac}(C/\langle \sigma\rangle)$ we can assume that $H_1$ and $H_2$ are not elliptic curves. 
From Hurwitz's genus formula and the fact that  $\sigma$ and $\iota\sigma$ 
are non-hyperelliptic involutions it follows that genus of the 
curve $H_1$ and $H_2$ is $2$. 
By easy Galois theory it follows that $C$ is the fibered product of
the two genus $2$ covers of $\PP^1$. Both $H_i$ are isomorphic
to the optimal curve given in previous subsection. If we represent them 
as double covers of $\PP^1$ five of the six ramification points must
coincide. But we know that in our case the ramification points
form two orbits of three points under the action of Frobenius. 
Therefore we have three or six common ramification points. This contradicts
our assumption.

Finally, we assume that the curve $C$ is not hyperelliptic, 
and therefore 
${\rm  Aut}(C) \cong {\rm  Aut}({\rm  Jac}(C))/\langle-1\rangle$ 
by the Theorem of Torelli. 
In case that the polarization does not corresponds to the modules of the cases $(2)$, $(3)$ or $(4)$ 
we have two commuting involutions $\sigma, \tau \in 
{\rm  Aut}( {\rm  Jac}(C), \Theta)/\langle-1\rangle$,
hence $\ZZ/2\ZZ \times \ZZ/2\ZZ$ in the automorphism group of $C$.
Thus we again obtain a diagram of the coverings:
$$
\xymatrix{
 &\ar[ld]_{2:1}C\ar[d]^{2:1}\ar[rd]^{2:1}&\\
 C/\langle\sigma\rangle\ar[rd]^{2:1}& C/\langle \sigma\tau\rangle \ar[d]^{2:1} & C/\langle \tau\rangle \ar[ld]_{2:1} \\
&C/\langle\sigma, \tau\rangle&\\}
$$
\par
If any quotient curve is not an optimal elliptic curve 
then each quotient curves $C/\langle \sigma \rangle$, 
$C/\langle \sigma  \tau \rangle$ and $C/\langle \tau \rangle$ has genus $2$.
Furthermore, $C/\langle \sigma,\; \tau\rangle$ cannot be $\PP^1$ 
since then $C$ is the fibered product of two optimal genus $2$ curves,
and this contradicts to the splitting of the Jacobian ${\rm  Jac}(C)$ as we saw above. So
$C/\langle \sigma,\tau\rangle$ must be an optimal elliptic curve $E$.

\par
We may assume that an optimal elliptic curve $E=C/\langle \sigma,\tau \rangle$  is given
by an equation $y^{2}=ax^{3}+bx^{2}+cx+d$ such that the optimal genus two curve  $C/\langle \sigma \rangle$ is given by equations
$$
\left\{
\begin{array}{l}
y^{2}=ax^{3}+bx^{2}+cx+d,\\
z^{2}=x.\\
\end{array}
\right.
$$

Now we show that any curve $H$ of genus $2$ having degree two map to  $E$ can be given by a system of equations
 $$
\left\{
\begin{array}{l}
y^{2}=ax^{3}+bx^{2}+cx+d,\\
z^{2}=\alpha_{0}+\alpha_{1}x+\alpha_{2}x^{2}+\beta y.\\
\end{array}
\right.
$$
where $\beta$ is either $1$ or some fixed element  from $\FF_q \setminus (\FF_{q})^{2}$.
Abusing notation we denote by $\infty$  a unique place of $E$ lying over $\infty \in \PP^{1}$ under $2:1$ map given by $x$ and let $D=f^{-1}(\infty)$ be a divisor of ${\rm Div}(H)$, where $f:H \rightarrow E$ is our degree $2$ map. 

If $\infty \in E$ is ramified then there exits a unique place $\infty^{\prime}$  of $\FF_{q}(H)$ lying above  the $\infty \in E$. Moreover, the Rimann-Roch vector space associated with the divisor $n\cdot \infty^{\prime}$ allows the following description in terms of basis:

\begin{tabular}{cl}
${\rm dim}(3\infty^{\prime})=2$ & $\{1,\, z \}$ with $(z)_{\infty}=2\infty^{\prime}$ 
(and hence $z \in \FF_{q}(H) \setminus \FF_{q}(E)$),\\
${\rm dim}(4\infty^{\prime})=3$ & $\{1,\, z, \, x\}$ with $(x)_{\infty}=3\infty^{\prime}$, \\
${\rm dim}(5\infty^{\prime})=4$ & $\{1,\, z, \, x, \, \tilde{z}\}$ with $(\tilde{z})_{\infty}=4\infty^{\prime}$ \\
\end{tabular}  

Due to the fact the vector space ${\rm L}(5\infty^{\prime})$ has dimension $5$, 
the elements $\{1,\, z,\, x,\, \tilde{x},\, z^2, y \}$ are linearly dependent over $\FF_{q}$. Thus either 
$\{1,\, z,\, x,\, z^2, y \}$ are linearly dependent (and hence $z^2=\beta_0+\beta_1x+\gamma y$) or
$\tilde{z}=\alpha_0+\alpha_1 x+\alpha_2 z +\alpha_3 y+\alpha_4 z^2$ with $\alpha_4,\alpha_3 \not =0.$  Taking into the account this expression, the facts that ${\rm dim\,}{\rm L}(8\infty^{\prime})=7$ and   
$\{1,\, z,\, x,\, \tilde{x},\,y, \,  z^2,\,  zx,\, x^2,\, \tilde{z}z\} 
\subset {\rm L}(8\infty^{\prime})$ it follows that $z$ is a root of the following polynomial
$$
 T^3+ \gamma_2 T^2+(\gamma_3 +\gamma_4y +
\gamma_5 x)T +\gamma_6 y+ \gamma_7 x^2+\gamma_8 x +\gamma_9=0,
$$
with $\gamma_4 \not = 0$. It yields that $z$ is an integral element over $\FF_q[x,y]$.
The degree of a minimal polynomial of $z$ is two and it divides the polynomial above. Hence
$$
 z^2+\delta_1 z +\delta_2 zy +
\delta_3 zx +\delta_4 y+ \delta_5 x^2+\delta_6 x +\delta_7=0,
$$
and the strict triangle inequality provides that $\delta_2=\delta_3=\delta_7=0$. As result we prove that there exists an element $z \in \FF_q(H)\setminus \FF_q(E)$ such that 
$z^2=\alpha_0+\alpha_1 x+\beta y.$

If the $\infty$ is unramified then $D=P_1+P_2$ over $\bar{\FF}_q.$ Using the similar technique for
$(z)_{\infty}=2P_1+P_2$, $(\tilde{z})_{\infty}=3P_1+2P_2$ and working with vector spaces
${\rm L}(3P_1+3P_2),\, {\rm L}(5P_1+4P_2)$, one obtains that
$$
z^2+\beta_1 zx +\beta_2 y +\beta_3 x^2+ \beta_4 x +\beta_5=0.
$$ 
As result there exists an element in $\FF_q(H)\setminus \FF_q(E)$, which is again denoted by $z$ such that
$z^2= \alpha_0+ \alpha_1 x +\alpha_2 x^2+\beta y.$

Now, we study the brach locus of $C/\langle  \tau\rangle \rightarrow E$ and $C/\langle \sigma \tau\rangle \rightarrow E$ on base that $C/\langle  \tau\rangle$ is given by equations
 $$
\left\{
\begin{array}{l}
y^{2}=ax^{3}+bx^{2}+cx+d,\\
z^{2}=\alpha_{0}+\alpha_{1}x+\alpha_{2}x^{2}+\beta y.\\
\end{array}
\right.
$$
If a place $P$ of $\FF_q(E)$ is ramified in the cover 
$C/\langle  \tau\rangle \rightarrow C/\langle \sigma, \tau\rangle $ then  $2$ does not divide $v_{P}(\alpha_{0}+\alpha_{1}x+\alpha_{2}x^{2}+\beta y)$, by the Kummer theory. On the other hand the curve $C/\langle \sigma \tau \rangle$ is given by
 $$
\left\{
\begin{array}{l}
y^{2}=ax^{3}+bx^{2}+cx+d,\\
z^{2}=x(\alpha_{0}+\alpha_{1}x+\alpha_{2}x^{2}+\beta y)\\
\end{array}
\right.
$$
 and $2$ divides $v_{P}(\alpha_{0}+\alpha_{1}x+\alpha_{2}x^{2}+\beta y)$ unless $P$ is a zero $x$
 (the zero of $x$, denoted by $P_0 \in \PP^1$, is unramified in the cover $E \rightarrow \PP^1$ given by $x$). So one of ramification point of $C\langle \tau \rangle \rightarrow E$ is a zero of $x$, otherwise the genus of $C\langle \tau \rangle$ is greater or equal to $3$ (since under this condition the cover $C\langle \sigma\tau \rangle \rightarrow E$ has al least three ramification points).
 
 Now we impose some conditions on the coefficients $\alpha_i$ via computing 
 $v_{Q_0}(\alpha_{0}+\alpha_{1}x+\alpha_{2}x^{2}+\beta y)$, where $Q_0$ is a place 
 of $\FF_q(C/\langle \sigma, \tau \rangle )=\FF_q(E)$ lying above $P_0 \in \PP^1$ (the zero of $x$).
 Note that one of the valuations  $v_{Q_0}((\alpha_{0}+\alpha_{1}x+\alpha_{2}x^{2}+\beta y))$ or
 $v_{Q_0}((\alpha_{0}+\alpha_{1}x+\alpha_{2}x^{2}-\beta y))$ is positive and hence
 $$
 \begin{array}{l}
 v_{Q_0}((\alpha_{0}+\alpha_{1}x+\alpha_{2}x^{2}+\beta y)(\alpha_{0}+\alpha_{1}x+\alpha_{2}x^{2}-\beta y))=\\
 v_{Q_0}((\alpha_{0}+\alpha_{1}x+\alpha_{2}x^{2})^2-\beta^2 y^2)=\\
 v_{Q_0}((\alpha_0^2-\beta^2d)+(2\alpha_0\alpha_1-\beta c )x+(\cdots)x^2 ) >0. 
 \end{array}
 $$
 Therefore one can get that $\alpha_0^2-\beta^2d=0$ (owing to the strict triangle inequality).

So for given an optimal elliptic curve $E$ with the prescribed conditions above 
we can run a computer search for the curve $C/\langle  \tau \rangle$ and $C/\langle \sigma \tau \rangle$ of the described form with
$\beta \in \{1,\, \mbox{a fixed non-square element of }\, \FF_q\}$ and two parameters $\alpha_1, \, \alpha_2 \in \FF_q$.

In additionally, we would like to remark  the property that for given an optimal elliptic curve
$y^2= \alpha x^3+\beta x^2+\gamma x+\delta$ there are at most two optimal elliptic curves of the form $y^2= (\alpha x^3+\beta x^2+\gamma x+\delta)(ux+w)$  implies that
there  are at most two candidates for the optimal elliptic curve  given by an equation $y^2=ax^3+bx^2+cx+d$ such that the elliptic curve  $y^2=(ax^3+bx^2+cx+d)x$ is also optimal.

A computer search shows that there are no such curves $C/\langle  \tau \rangle$, $C/\langle \sigma \tau \rangle$  and hence there is no an optimal curve of genus $4$ satisfying the diamond diagram above, if $q<10^3$.

\par
If one of the quotient curves is an optimal elliptic curve then from the splitting of the Jacobian of curve $C$ 
for the Klein $4$-group 
it follows that the second one is also an optimal elliptic curve and the third quotient curve
is an optimal curve of genus $2$; moreover $C/\langle \sigma, \tau\rangle$ is isomorphic to $\PP^1$.
As a result we have that an optimal curve of genus $4$ is the fibered product of optimal elliptic curves if and 
only if an optimal curve of genus $2$ can be given by an equation $z^2=F(x)$ such that $F(x)$ splits into two irreducible polynomials $g_1(x)$, $g_2(x)$ of degree $3$ and such that the equations $u^2=g_1(x)$, $w^2=g_2(x)$ give optimal elliptic curves.
An easy computer check for each $q$ shows that an optimal curve of genus $2$ is not given by such an equation. 
\end{proof}
\end{subsection}
\begin{subsection}{Optimal curves of genus $5$.}
\begin{lemma}\label{degree}
An optimal curve $C$ of genus $5$ over $\FF_{q}$ with discriminant $-19$ is not a 
hyperelliptic curve nor a double cover
of a curve of genus $\geq 3$. 
\end{lemma} 
\begin{proof}
If $C$ is a double cover of a genus $3$ curve
then it is an unramified cover and its number of rational points is even. But $m=[2\sqrt{q}]$ is odd and
hence $\#C(\FF_{q})=q+1\pm 5\cdot m$ is odd.
The Hurwitz formula excludes the case that it is a double cover of a curve of genus $\geq 4$.
\par
If $C$ is a hyperelliptic curve with the hyperelliptic involution $\tau$ then from the structure of the automorphism group of the irreducible unimodular hermitian modules it follows that there is an extra involution. Let us denote this involution by $\sigma$.  Then we have the following relation on Jacobians 
$E^5\cong {\rm  Jac}(C)\sim {\rm  Jac}(C/\langle\sigma\rangle)\times {\rm  Jac}(C/\langle-\sigma\rangle)$,
where $E$ is an optimal elliptic curve.
Therefore one of the quotient curves is a curve of genus either $4$ or $3$ and this is impossible. 
\end{proof}
Now we show that the automorphism group of a maximal curve of 
genus~$5$ over $\FF_{q}$ cannot have  the Klein $4$-group $\ZZ/2\ZZ\times \ZZ/2\ZZ$ as subgroup.
\begin{lemma}\label{order4}
If $C$ is a maximal optimal curve of genus $5$ over $\FF_{q}$ with discriminant $-19$ 
then ${\rm   Aut}(C)$ does not contain  the Klein $4$-group.
\end{lemma}
\begin{proof}
Suppose that $\sigma$ and $\tau$ are two commuting
involutions of $C$. Then $C$ is the fibered product of
$C_1:=C/\langle \sigma \rangle$ and $C_2:=C/\langle \tau \rangle$
over $C/\langle \sigma, \tau \rangle$. Put $C_3:= C/\langle \sigma
\tau \rangle$. Then we have the isogeny
$$
{\rm   Jac}(C_1)\times {\rm   Jac}(C_2) \times {\rm   Jac}(C_3)\sim
{\rm   Jac}(C/\langle \sigma, \tau \rangle)^2 \times{\rm   Jac}(C).
$$
In view of Proposition~\ref{degree} the dimension of the 
part on the left-hand-side is $\leq 6$, but the dimension of the right-hand-side
is odd and $\geq 5$. The only possibility is that
one of the $C_i$, say $C_1$, has genus~$1$ and the other two have genus $2$
and that $C/\langle \sigma, \tau \rangle$ has genus $0$.
But then the $C_2$ and $C_3$ share two of their six ramification points
and this contradicts  Remark~\ref{weierstrass2},
which implies that the six ramification points decompose in two Galois orbits of length~$3$.
\end{proof}
\begin{theorem}
If there exits an optimal curve $C$ of genus $5$ over a finite field $\FF_{q}$ with
$d(\FF_{q})=-19$ then its automorphism group is isomorphic to the Dihedral group  ${\rm D}_{5}$.
\end{theorem}
\begin{proof}
Looking at generators of automorphism groups of unimodular irreducible
hermitian modules, we see that ${\rm  Aut}_{\FF_{q}}(C)$ contains either the Klein $4$-group,
or it has order $10$ (cf. the Appendix). Due to the Lemma~\ref{order4} 
we may exclude all automophism groups with $2$-Sylow subgroup of order greater or equal to $2^{3}$ since it has  the Klein $4$-group as a subgroup.  On the other hand neither
an optimal elliptic curve nor an optimal curve of genus $2$ has no an auromorphism of order $5$ and hence ${\rm  Aut}_{\FF_{q}}(C) \cong {\rm D}_{5}$ . 
\end{proof}
\end{subsection}
\end{section}

\begin{section}{Appendix}
Here we present the unimodular irreducible hermitian modules of dimension $4$ and $5$ 
over an imaginary quadratic extension $K$ of $\QQ$ with discriminant $d(K)=-19$ such that
their automorphism group does not contain $\ZZ/5\ZZ$ as a subgroup. We also show that
the Klein $4$-group is a subgroup of each automorphism group.
The irreducible hermitian modules and elements the automorphism groups were taken from
\cite{Schulze-letter}.
\par
Throughout we let $w=\frac{1+\sqrt{-19}}{2}$.
\begin{subsection}{Dimension $4$.}
\footnotesize
\begin{center}
$
H_1=\left(
\begin{array}{cccc}
       3&  & &\\
      -1&        3& & \\
      -1&      1+w&        3& \\
     2-w&        0&       -1&        3\\
\end{array}
\right),
\;  
|{\rm Aut}| = 2^4\cdot3,
$
\end{center}
$
\alpha_2=\left(
\begin{array}{cccc}
       0&        0&        1&        0\\
       0&      1+w&       -3&       -1\\
      -1&     -2+w&     -1-w&        0\\
       0&        1&        0&        0\\
\end{array}
\right),
$
$
\alpha_3=\left(
\begin{array}{cccc}
       1&        0&        0&        0\\
      -3&       -1&        0&      1+w\\
    -1-w&        0&       -1&     -2+w\\
       0&        0&        0&        1\\
\end{array}
\right).
$
$$
\alpha_2^4=-1,\; \alpha_3^2=1,\; (\alpha_2)^2\cdot \alpha_3 =-\alpha_3 \cdot (\alpha_2)^2.
$$

\begin{center}
$
H_5=\left(
\begin{array}{cccc}
       2& & &\\
       0&        2& & \\
       0&       -1&        2&\\
      -1&     -1+w&      1-w&        4\\ 
 \end{array}
\right),
\;     
|{\rm Aut}| = 2^3 \cdot 3^2$,
\end{center}
$
\alpha_3=\left(
\begin{array}{cccc}
      -2&     -1+w&      1-w&       -3\\
       0&        1&        0&        0\\
       0&       -1&       -1&        0\\
       1&        0&     -1+w&        2\\
\end{array}
\right),
$
$
\alpha_5=\left(
\begin{array}{cccc}
      -1&     -1+w&      1-w&       -3\\
       0&        1&        0&        0\\
       0&        0&        1&        0\\
       0&        0&        0&        1\\
\end{array}
\right).
$
$$
(\alpha_3  \cdot \alpha_5)^6=\alpha_5^2=1,\; (\alpha_3  \cdot \alpha_5)^3 \cdot \alpha_5 =\alpha_5 \cdot (\alpha_3  \cdot \alpha_5)^3.
$$

\begin{center}
$
H_6=\left(
\begin{array}{cccc}
       2& & &\\
      -1&        2& & \\
       0&       -1&        3& \\
       0&       -w&     -1+w&        4 \\
 \end{array}
\right),
\;  
|{\rm Aut}| = 2^3\cdot 3$,
\end{center}
$
\alpha_2=\left(
\begin{array}{cccc}
       0&        1&        0&        0\\
      -1&       -1&        0&        0\\
     2-w&     3-2w&       -1&     -2-w\\
       w&        w&        0&        1\\
\end{array}
\right),
$
$
\alpha_3=\left(
\begin{array}{cccc}
       0&        1&        0&        0\\
       1&        0&        0&        0\\
    -2+w&    -3+2w&        1&      2+w\\
      -w&       -w&        0&       -1\\
\end{array}
\right).
$
$$
\alpha_2^6=\alpha_3^2=1,\; \alpha_3 \cdot \alpha_2^3=\alpha_2^3 \cdot \alpha_3.
$$

\begin{center}
$
H_7=\left(
\begin{array}{cccc}
       2& & &\\
       1&        2& & \\
      -1&       -1&        3& \\
    -1+w&     -1+w&       -w&        4\\
 \end{array}
\right),
\;
|{\rm Aut}| = 2^3\cdot3,
$
\end{center}
$
\alpha_3=\left(
\begin{array}{cccc}
       1&        0&        0&        0\\
       1&       -1&        0&        0\\
      -1&        0&       -1&        0\\
    -1+w&        0&        0&       -1\\
\end{array}
\right),
$
$
\alpha_4=\left(
\begin{array}{cccc}
       1&        0&        0&        0\\
       0&        1&        0&        0\\
       w&        w&       -1&     -3+w\\
       0&        0&        0&        1\\
\end{array}
\right).
$
$$
\alpha_3^2=\alpha_4^2=1, \; \alpha_3 \cdot \alpha_4=\alpha_4 \cdot \alpha_3.
$$

\begin{center}
$
H_8=\left(
\begin{array}{cccc}
       2& & &\\
       0&        2& & \\
      -w&        0&        3& \\
       0&       -w&        0&        3\\
 \end{array}
\right),
\;
|{\rm Aut}| = 2^5,
$
\end{center}
$
\alpha_2=\left(
\begin{array}{cccc}
       1&        0&        0&        0\\
       0&       -1&        0&        0\\
      -w&        0&       -1&        0\\
       0&        0&        0&       -1\\
\end{array}
\right),
$
$
\alpha_3=\left(
\begin{array}{cccc}
       1&        0&        0&        0\\
       0&       -1&        0&        0\\
       0&        0&        1&        0\\
       0&        0&        0&       -1\\
\end{array}
\right).
$
$$
\alpha_2^2=\alpha_2^2=1, \; \alpha_2 \cdot \alpha_3=\alpha_3 \cdot \alpha_2.
$$

\begin{center}
$
H_9=\left(
\begin{array}{cccc}
       2& & &\\
       0&        2& & \\
       0&       -1&        3& \\
      -1&       -1&      1-w&        3\\
\end{array}
\right),
\;
|{\rm Aut}| = 2^4,
$
\end{center}
$
\alpha_3=\left(
\begin{array}{cccc}
      -1&        0&        0&        0\\
       0&        1&        0&        0\\
       0&        0&        1&        0\\
       1&        0&        0&        1\\
\end{array}
\right),
$
$
\alpha_4=\left(
\begin{array}{cccc}
       1&        0&        0&        0\\
       0&        1&        0&        0\\
       0&       -1&       -1&        0\\
      -1&       -1&        0&       -1\\
\end{array}
\right).
$

$$
\alpha_3^2=\alpha_4^2=1, \; \alpha_3 \cdot \alpha_4=\alpha_4 \cdot \alpha_3.
$$

\end{subsection}
\begin{subsection}{Dimension $5$.}
\footnotesize
\begin{center}
$
H_1=\left(
\begin{array}{ccccc}
       3&       &       &        &        \\
      -1&        3&        && \\     
      -1&        1&        3&        &  \\
       0&       -w&        0&        3&  \\
       1&     -1+w&     -2+w&       -1&        4\\
 \end{array}
\right),
\;
|{\rm Aut}| = 2^3,
$
\end{center}
$
\alpha_2=\left(
\begin{array}{ccccc} 
    -1-w&     -2w&      2-w&       -2&        4\\
       0&       1&        0&        0&        0\\
       w&      2w&     -3+w&        2&       -4\\
       0&       0&        0&        1&        0\\
      -w&      -w&      2-w&       -1&        3\\
\end{array}
\right),
$\\
$
\alpha_3=\left(
\begin{array}{ccccc}
       0&        0&        1&        0&        0\\
      -w&     1-2w&      2-w&       -2&        4\\
       1&        0&        0&        0&        0\\
       0&        0&        0&        1&        0\\
       2&        5&      2+w&      1-w&    -1+2w\\
 \end{array}
\right).
$
$$
\alpha_2^2=\alpha_3^2=1, \; \alpha_2\cdot \alpha_3=\alpha_3\cdot \alpha_2.
$$

\begin{center}
$              
H_6=\left(
\begin{array}{ccccc}
       2& & & &\\
       1&        2& &&\\
      -1&       -1&        3& &\\
       1&        1&        0&        3& \\
       0&        0&        w&       -1&        3\\
\end{array}
\right),
\;
|{\rm Aut}| = 2^4 \cdot 3,
$
\end{center}
$              
\alpha_1=\left(
\begin{array}{ccccc}
       1&       -1&        0&        0&        0\\
       1&        0&        0&        0&        0\\
      -1&        0&       -1&        0&        0\\
       1&        0&        0&       -1&        0\\
       0&        0&        0&        0&       -1\\
\end{array}
\right),
$
$              
\alpha_2=\left(
\begin{array}{ccccc}
       w&        w&    -2+2w&     -1-w&     -3-w\\
       w&     -1+w&    -2+2w&     -1-w&     -3-w\\
      -w&       -w&     1-2w&      1+w&      3+w\\
       1&        0&        0&       -1&        0\\
       0&        0&        0&        0&       -1\\
\end{array}
\right),
$
$              
\alpha_3=\left(
\begin{array}{ccccc}
       w&        w&    -2+2w&     -1-w&     -3-w\\
       0&        1&        0&        0&        0\\
       0&        0&        1&        0&        0\\
    -1+w&        w&    -2+2w&       -w&     -3-w\\
       0&        0&        0&        0&        1\\
\end{array}
\right).
$
$$
(\alpha_2\cdot \alpha_3)^2=(\alpha_1 \cdot \alpha_3)^4=1, \; 
(\alpha_2\cdot \alpha_3)\cdot (\alpha_1 \cdot \alpha_3)^2=(\alpha_1 \cdot \alpha_3)^2\cdot (\alpha_2\cdot \alpha_3)
$$

\begin{center}
$              
H_7=\left(
\begin{array}{ccccc}
       2& & & &\\
       1&        2& &&\\
      -1&        0&        3& &\\
       1&        0&        0&        3& \\
       0&        0&      1-w&       -1&        3\\ 
\end{array}
\right),
\;
|{\rm Aut}| = 2^4\cdot 3,
$
\end{center}
$             
\alpha_2=\left(
\begin{array}{ccccc}
      -1&        1&        0&        0&        0\\
   -2+2w&      1-w&       2w&      2-w&      4-w\\
       1&       -1&        1&        0&        0\\
     -2w&        w&      -2w&     -1+w&     -4+w\\
       0&        0&        0&        0&        1\\
\end{array}
\right),
$\\
$              
\alpha_3=\left(
\begin{array}{ccccc}
   -1+2w&       -w&       2w&      2-w&      4-w\\
   -2+2w&      1-w&       2w&      2-w&      4-w\\
      -1&        1&       -1&        0&        0\\
      2w&       -w&       2w&      1-w&      4-w\\
       0&        0&        0&        0&       -1\\
\end{array}
\right).
$
$$
(\alpha_2)^4=\alpha_3^2=1, \; (\alpha_2)^2\cdot \alpha_3=\alpha_3 \cdot (\alpha_2)^2.
$$

\begin{center}
$
H_8=\left(
\begin{array}{ccccc}
       2& & & &\\
       0&        2& & &\\
       0&       -1&        3& &\\
       0&        1&     -1-w&        3& \\
     1-w&        0&        0&        0&        3\\
\end{array}
\right),
\;
|{\rm Aut}| = 2^4 \cdot 3,
$
\end{center}
$              
\alpha_2=\left(
\begin{array}{ccccc}
       1&        0&        0&        0&        0\\
       0&       -2&     -2+w&      1+w&        0\\
       0&        1&      1-w&     -1-w&        0\\
       0&       -1&     -2+w&        w&        0\\
     1-w&        0&        0&        0&       -1\\
\end{array}
\right),
$
$              
\alpha_3=\left(
\begin{array}{ccccc}
       1&        0&        0&        0&        0\\
       0&       -2&     -2+w&      1+w&        0\\
       0&        1&      1-w&     -1-w&        0\\
       0&       -1&     -2+w&        w&        0\\
       0&        0&        0&        0&        1\\
\end{array}
\right).
$
$$
\alpha_2^2=-1, \; \alpha_3^2=1, \; \alpha_2 \cdot \alpha_3 =\alpha_3 \cdot \alpha_2.
$$

\begin{center}
$
H_9=\left(
\begin{array}{ccccc}
       2& & & &\\
      -1&        2& && \\
       0&       -1&        3& &\\
      -1&        1&        0&        3& \\
      -1&        0&        1&      1-w&        3\\
\end{array}
\right),
\;
|{\rm Aut}| = 2^3\cdot 3,
$
\end{center}
$              
\alpha_3=\left(
\begin{array}{ccccc}
      -1&       -1&        0&        0&        0\\
       0&        1&        0&        0&        0\\
      -3&     -1+w&        1&    -1-2w&     -5+w\\
       0&        1&        0&       -1&        0\\
       0&        0&        0&        0&       -1\\
\end{array}
\right),
$
$              
\alpha_4=\left(
\begin{array}{ccccc}
       1&        0&        0&        0&        0\\
       0&        1&        0&        0&        0\\
       3&       -w&       -1&     1+2w&      5-w\\
       0&        0&        0&        1&        0\\
       0&        0&        0&        0&        1\\
\end{array}
\right).
$
$$
\alpha_3^2=\alpha_4^2=1, \; \alpha_3\cdot \alpha_4= \alpha_4 \cdot \alpha_3.
$$

\begin{center}
$
H_{10}=\left(
\begin{array}{ccccc}
       2& & & &\\
       1&        2& & & \\
       1&        0&        3& & \\
       0&        1&       -w&        3& \\
       1&        0&        1&       -1&        3\\
\end{array}
\right),
\;
|{\rm Aut}| = 2^3 \cdot 3,$
\end{center}
$
\alpha_3=\left(
\begin{array}{ccccc}
       0&       -1&        0&        0&        0\\
      -1&        0&        0&        0&        0\\
       1&       -1&       -1&        0&        0\\
      -1&        1&        0&       -1&        0\\
     3-w&       -4&    -3+2w&      4+w&        1\\
\end{array}
\right),
$ 
$              
\alpha_4=\left(
\begin{array}{ccccc}
       1&        0&        0&        0&        0\\
       0&        1&        0&        0&        0\\
       0&        0&        1&        0&        0\\
       0&        0&        0&        1&        0\\
    -2+w&        3&     3-2w&     -4-w&       -1\\
\end{array}
\right).
$
$$
\alpha_3^2=\alpha_4^2=1, \; \alpha_3\cdot \alpha_4= \alpha_4 \cdot \alpha_3.
$$

\begin{center}
$
H_{11}=\left(
\begin{array}{ccccc}
       2& & & &\\
       0&        3& & &\\
       0&        1&        3& & \\
       0&      2-w&        0&        3& \\
      -1&        0&        1&       -1&        3\\
\end{array}
\right),
\;
|{\rm Aut}| = 2^3\cdot 3,$
\end{center}
$              
\alpha_2=\left(
\begin{array}{ccccc}
       1&    -3+3w&       -w&      7-w&        3\\
       0&        1&        0&        0&        0\\
       0&        3&       -1&     -1-w&        0\\
       0&        0&        0&        1&        0\\
       0&      2-w&        0&       -3&       -1\\
\end{array}
\right),
$ 
$              
\alpha_3=\left(
\begin{array}{ccccc}
       1&        0&        0&        0&        0\\
       0&       -1&        0&        0&        0\\
       0&       -3&        1&      1+w&        0\\
       0&        0&        0&       -1&        0\\
       0&     -2+w&        0&        3&        1\\
\end{array}
\right).
$
$$
\alpha_2^2=\alpha_3^2=1, \; \alpha_3\cdot \alpha_2= \alpha_2 \cdot \alpha_3.
$$

\begin{center}
$
H_{12}=\left(
\begin{array}{ccccc}
       2& & & &\\
       0&        3& & &\\
       0&        1&        3& &\\
       0&       -1&     -1-w&        3& \\
       1&       -w&       -1&        0&        3 \\
\end{array}
\right),
\;
|{\rm Aut}| = 2^3 \cdot 3,
$
\end{center}
$              
\alpha_2=\left(
\begin{array}{ccccc}
       1&        0&        0&        0&        0\\
       0&       -1&      2-w&     -1-w&        0\\
       0&        0&        1&        0&        0\\
       0&        0&        0&        1&        0\\
       1&        0&       -3&     -2+w&       -1\\
\end{array}
\right),
$ 
$              
\alpha_3=\left(
\begin{array}{ccccc}
       1&        0&        0&        0&        0\\
       0&        1&     -2+w&      1+w&        0\\
       0&        0&       -1&        0&        0\\
       0&        0&        0&       -1&        0\\
       0&        0&        3&      2-w&        1\\
\end{array}
\right).
$ 
$$
\alpha_2^2=\alpha_3^2=1, \; \alpha_3\cdot \alpha_2= \alpha_2 \cdot \alpha_3.
$$

\begin{center}
$
H_{13}=\left(
\begin{array}{ccccc}
       2& & & & \\
       0&        3& & &\\
      -1&       -w&        3& &\\
      -1&       -1&        0&        3& \\
     1-w&        1&     -1+w&        0&        4\\
\end{array}
\right),
\;
|{\rm Aut}| = 2^3,$
\end{center}
$             
\alpha_2=\left(
\begin{array}{ccccc}
       1&        0&        0&        0&        0\\
   -1+2w&     -2+w&      1+w&     -1+w&        2\\
       0&        0&        1&        0&        0\\
    1-2w&      3-w&     -1-w&      2-w&       -2\\
       0&        0&        0&        0&        1\\
\end{array}
\right),
$\\ 
$              
\alpha_3=\left(
\begin{array}{ccccc}
       1&        0&        0&        0&        0\\
    1-2w&      2-w&     -1-w&      1-w&       -2\\
      -1&        0&       -1&        0&        0\\
   -2+2w&     -3+w&      1+w&     -2+w&        2\\
     1-w&        0&        0&        0&       -1\\
\end{array}
\right).
$
$$
\alpha_2^2=\alpha_3^2=1, \; \alpha_3\cdot \alpha_2= \alpha_2 \cdot \alpha_3.
$$

\begin{center}
$
H_{14}=\left(
\begin{array}{ccccc}
       2& & & &\\
      -1&        3& & &\\
       0&       -1&        3& &\\
       0&      1-w&       -1&        3& \\
       0&       -1&      1-w&        1&        3\\
\end{array}
\right),
\;
|{\rm Aut}| = 2^4,
$
\end{center}
$              
\alpha_2=\left(
\begin{array}{ccccc}
      -1&       -2&        2&      1+w&     -1-w\\
       0&        0&       -1&        0&        0\\
       0&       -1&        0&        0&        0\\
       0&        0&        0&        0&       -1\\
       0&        0&        0&       -1&        0\\
\end{array}
\right),
$ 
$              
\alpha_4=\left(
\begin{array}{ccccc}
       1&        0&        0&        0&        0\\
      -1&       -1&        0&        0&        0\\
      -1&       -2&        1&      1+w&     -1-w\\
       0&        0&        0&       -1&        0\\
       0&        0&        0&        0&       -1\\
\end{array}
\right).
$
$$
\alpha_2^2=\alpha_4^2=1, \; \alpha_4\cdot \alpha_2= \alpha_2 \cdot \alpha_4.
$$

\begin{center}
$
H_{15}=\left(
\begin{array}{ccccc}
       2& & & & \\
       0&        2& & & \\
       1&       -1&        3& & \\
       0&        0&     -1+w&        3& \\
       w&        0&        w&        1&        4\\
\end{array}
\right),
\; 
|{\rm Aut}| = 2^4,$
\end{center}
$              
\alpha_2=\left(
\begin{array}{ccccc}
       0&        1&        0&        0&        0\\
       1&        0&        0&        0&        0\\
       0&        0&       -1&        0&        0\\
       0&        0&        0&       -1&        0\\
       1&       -1&     -2-w&      1-w&        1\\
\end{array}
\right),
$ 
$              
\alpha_4=\left(
\begin{array}{ccccc}
       1&        0&        0&        0&        0\\
       0&        1&        0&        0&        0\\
       1&       -1&       -1&        0&        0\\
       0&        0&        0&       -1&        0\\
       w&        0&        0&        0&       -1\\
\end{array}
\right).
$
$$
\alpha_2^2=\alpha_4^2=1, \; \alpha_4\cdot \alpha_2= \alpha_2 \cdot \alpha_4.
$$

\begin{center}
$
H_{16}=\left(
\begin{array}{ccccc}
       2& & & &\\
       0&        2& & &\\
      -1&        1&        3& &\\
       0&        0&       -w&        3& \\
       0&     -1+w&     -1+w&       -1&        4\\ 
\end{array}
\right),
\;
|{\rm Aut}| = 2^4,$
\end{center}
$              
\alpha_2=\left(
\begin{array}{ccccc}
       0&        1&        0&        0&        0\\
       1&        0&        0&        0&        0\\
       1&       -1&        1&        0&        0\\
       0&        0&        0&        1&        0\\
    -2+w&        1&     -3+w&        w&       -1\\
\end{array}
\right),
$ 
$              
\alpha_4=\left(
\begin{array}{ccccc}
       1&        0&        0&        0&        0\\
       0&        1&        0&        0&        0\\
      -1&        1&       -1&        0&        0\\
       0&        0&        0&       -1&        0\\
       0&     -1+w&        0&        0&       -1\\
\end{array}
\right).
$ 
$$
\alpha_2^2=\alpha_4^2=1, \; \alpha_4\cdot \alpha_2= \alpha_2 \cdot \alpha_4.
$$

\begin{center}
$
H_{17}=\left(
\begin{array}{ccccc}
       2& & &  &\\
      -1&        3& & &\\\
       0&       -1&        3& &\\
       0&     -1+w&       -1&        3& \\
       0&        0&       -1&        0&        3 \\
\end{array}
\right),
\;
|{\rm Aut}| = 2^4,$
\end{center}

$              
\alpha_2=\left(
\begin{array}{ccccc}
     4-w&     8-2w&        5&     5+2w&        2\\
    -4+w&    -8+2w&       -5&    -5-2w&       -1\\
       0&        0&        1&        0&        0\\
      -w&      -2w&      1-w&      4-w&        0\\
       1&        1&        0&        0&        0\\
\end{array}
\right),
$ 
$              
\alpha_4=\left(
\begin{array}{ccccc}
       1&        0&        0&        0&        0\\
      -1&       -1&        0&        0&        0\\
       0&        0&       -1&        0&        0\\
       0&        0&        0&       -1&        0\\
     4-w&     8-2w&        5&     5+2w&        1\\
\end{array}
\right).
$
$$
\alpha_2^2=\alpha_4^2=1, \; \alpha_4\cdot \alpha_2= \alpha_2 \cdot \alpha_4.
$$

\begin{center}
$
H_{18}=\left(
\begin{array}{ccccc}
       2& & & &\\
       0&        3& & &\\
      -1&       -1&        3& &\\
       0&      1-w&        w&        3& \\
       0&        1&        1&        1&        3 \\
\end{array}
\right),
\;
|{\rm Aut}| = 2^4,$
\end{center}
$              
\alpha_3=\left(
\begin{array}{ccccc}
      -1&        0&        0&        0&        0\\
       0&        1&        0&        0&        0\\
       1&        0&        1&        0&        0\\
       0&        0&        0&        1&        0\\
       0&        0&        0&        0&        1\\
\end{array}
\right),
$ 
$              
\alpha_4=\left(
\begin{array}{ccccc}
       1&        0&        0&        0&        0\\
       0&       -1&        0&        0&        0\\
      -1&        0&       -1&        0&        0\\
       0&        0&        0&       -1&        0\\
    -1+w&      2-w&    -2+2w&     -5-w&        1\\
\end{array}
\right).
$ 
$$
\alpha_3^2=\alpha_4^2=1, \; \alpha_4\cdot \alpha_3= \alpha_3 \cdot \alpha_4.
$$

\begin{center}
$
H_{19}=\left(
\begin{array}{ccccc}
       2& & & &\\
       0&        3& & &\\
       0&        0&        3& &\\
       0&       -w&       -1&        3& \\
       1&     -1+w&        0&       -1&        3 \\
\end{array}
\right),
\;
|{\rm Aut}| = 2^4,$
\end{center}

$              
\alpha_3=\left(
\begin{array}{ccccc}
      -1&        0&        0&        0&        0\\
       0&        1&        0&        0&        0\\
       0&        0&        1&        0&        0\\
       0&        0&        0&        1&        0\\
      -1&        0&        0&        0&        1\\
\end{array}
\right),
$
$              
\alpha_4=\left(
\begin{array}{ccccc}
       1&        0&        0&        0&        0\\
       0&       -1&        0&        0&        0\\
     2-w&     2+3w&        1&        5&    -4+2w\\
       0&        0&        0&       -1&        0\\
       1&        0&        0&        0&       -1\\
\end{array}
\right).
$ 
$$
\alpha_3^2=\alpha_4^2=1, \; \alpha_4\cdot \alpha_3= \alpha_3 \cdot \alpha_4.
$$

\begin{center}
$
H_{20}=\left(
\begin{array}{ccccc}
       3& & & &\\
       1&        3& & &\\
       0&     -2+w&        3& &\\
      -1&        0&        0&        3& \\
      -w&       -1&        1&       -1&        3 \\
\end{array}
\right),
\;
|{\rm Aut}| = 2^3,$
\end{center}

$              
\alpha_2=\left(
\begin{array}{ccccc}
      -1&        3&      1+w&        0&        0\\
       0&        1&        0&        0&        0\\
       0&        0&        1&        0&        0\\
     2+w&      1-w&        1&        1&      3-w\\
       0&     -1-w&      2-w&        0&       -1\\
\end{array}
\right),
$ 
$              
\alpha_3=\left(
\begin{array}{ccccc}
       1&       -3&     -1-w&        0&        0\\
       0&       -1&        0&        0&        0\\
       0&        0&       -1&        0&        0\\
       0&        0&        0&        1&        0\\
       0&      1+w&     -2+w&        0&        1\\
\end{array}
\right).
$
$$
\alpha_2^2=\alpha_3^2=1, \; \alpha_3\cdot \alpha_2= \alpha_2 \cdot \alpha_3.
$$

\end{subsection}

\end{section}


\begin{thebibliography}{99}
\bibitem{Andreotti} Aldo Andreotti: On a Theorem of Torelli.
{\sl American \ J.\ of Math.\ {\bf 80} no.\  4}, pp.\ 801-828 (1958).
\bibitem{Deligne} Pierre Deligne: Variet\'es ab\'eliennes ordinaires sur un corps fini.
{\sl Invent.\ Math.\  {\bf 8} }, pp.\ 238-243 (1969).
\bibitem{Geer} Gerard van der Geer and Marcel van der Vlugt: Supersingular curves of genus 2 over finite fields of characteristic 2, {\sl  Math.\  Nachr.\ {\bf 159} }, pp.\ 73-81 (1992).
\bibitem{Howe} Everett W.\ Howe: Principally polaried ordinary abelian varieties over finite field.
{\sl  Trans.\ Amer.\ Math.\ Soc.\ {\bf 347} No.\bf {7}}, pp.\ 2361-2401 (1995).
\bibitem{KR} E.\ Kani and M.\ Rosen: Idempotent relations and factors of Jacobians. 
{\sl  Math.\ Ann.\  {\bf 284} no.\ 2 }, pp.\ 307-327 (1989).
\bibitem{Lauter1} K.\ Lauter: The maximum or minimum number
of rational points on curves of genus three curves over finite fields.  (With 
an appendix by J.\-P.\ Serre). {\sl  Compositio Math.\ {\bf 134} }, pp.\ 87-111 (2002).
\bibitem{Lauter2} K.\ Lauter: Geometric methods for improving the upper bounds on the number of rational points on algebraic curves over finite fields.  (With 
an appendix by J.\-P.\ Serre). {\sl  J. Algebraic Geom.\ {\bf 10} }, pp.\ 19-36 (2001).
\bibitem{Oort} F.\ Oort and K.\ Ueno: Principally polarized abelian varieties of 
dimension two or three are Jacobian varieties.
{\sl  J.\ Fac.\ Sci.\ Univ.\ Tokyo {\bf 20} }, pp.\ 377-381 (1973).
\bibitem{Roquette} P.\ Roquette: 
Absch\"atzung der Automorphismenanzahl von Funktionenk\"orpern bei Primzahlcharakteristik,
{\sl  Math.\ Z.\ {\bf 117} }, pp.\ 157-163 (1970).
\bibitem{Sch} A.\ Schiemann: Classification of Hermitian Forms with the Neighbour Method. 
{\tt http://www.math.uni-sb.de/ag/schulze/Hermitian-lattices/}.
\bibitem{Scht} H.\ Stichtenoth: Algebraic Function Fields and Codes,
{\sl  Springer Universitext  Berlin-Heidelberg}, (1993).
\bibitem{Schulze-letter} R.\ Schulze-Pillot: Letter to A. Zaytsev, 3.\ 09.\ 2006.
\bibitem{Serre2} J-P.\ Serre: Nombre de points des courbes  alg\'ebriques sur $\FF_q$,
{\sl S\'em. de Th\'eorie des Nombres de Bordeaux} , exp.\ no.\ 22.\ 
(= Oeuvres III, No.\ 129, pp.\ 664-668), (1982/83).
\bibitem{Serre3} J-P.\ Serre: R\'esum\'e des Cours de 1983-1984, 
(= Oeuvres III, No.\ 132, pp.\ 701-705).
\bibitem{Serre}  Appendix by J.-P. Serre: in \cite{Lauter1}.
\bibitem{Shabat} V.\ Shabat, Ph.\ D.\ thesis: Curves with Many Points, Amsterdam 2001,
{\tt http://www.science.uva.nl/math/Research/Dissertations/Shabat2001.text.pdf}
\bibitem{Singh} B.\ Singh: On the group of automorphisms of a function field of genus at least two, 
{\sl  J.\ of Pure and Applied Algebra {\bf 4} }, pp.\ 203-229 (1974).
\bibitem{Nakajima}  S. Nakajima: $p$-rank and automorphism group of algebraic curves,
{\sl Transactions of the American Mathematical Society}, vol.\ 303, no.\ 2, pp.\ 595-607.


\bibitem{Waterhouse} William C.\ Waterhouse: Abelian varieties over finite fields,
{\sl  Ann.\ Sci.\ Ecole Norm.\ Sup.\ {\bf 2} }, pp.\ 521-560 (1969).
\end{thebibliography}
\end{document}